\documentclass[12pt,righttag]{amsart}
\usepackage{amsmath,amsthm,amssymb,amscd, amstext, amsopn, amsxtra}
\numberwithin{equation}{section}
\theoremstyle{definition}
  \newtheorem{theorem}{Theorem}[section]
  
  \newtheorem{proposition}[theorem]{Proposition}
  \newtheorem{lemma}[theorem]{Lemma}
  \newtheorem {corollary}[theorem]{Corollary}

\theoremstyle{remark}

\newcommand{\etalchar}[1]{$^{#1}$}

\newcommand{\I}{I}
\newcommand{\ee}[1]{\tilde{e}_{{#1}} }
\newcommand{\f}[1]{\tilde{f}_{{#1}} }
\newcommand{\ff}[1]{\tilde{f}_{{#1}} }
\newcommand{\ffi}{\f{i}}
\newcommand{\fj}{\f{j}}
\newcommand{\ei}{\ee i}
\newcommand{\ej}{\ee j}
\newcommand{\eee}[1]{{e}_{{#1}} }
\newcommand{\eei}{\eee i}
\newcommand{\eej}{\eee j}
\newcommand{\eps}{{\varepsilon} }
\newcommand{\epsof}[1]{{\varepsilon}_{{#1}} }
\newcommand{\epsi}{\epsof i}

\newcommand{\omitt}[1]{}
 
\newcommand{\aA}{{\mathbf a}}
\newcommand{\bb}{{\mathbf b}}
\newcommand{\cc}{{\mathbf c}}
\newcommand{\dd}{{\mathbf d}}
\newcommand{\zero}{{\mathbf 0}}
\newcommand{\singular}{singular}
\newcommand{\walk}{walk}
\newcommand{\walks}{walks}
\newcommand{\wt}{{\mathrm wt}}
\newcommand{\fin}{{\mathrm fin}}

\newcommand{\Z}{{\mathbb Z}}

\DeclareMathOperator{\Rep}{Rep}
\DeclareMathOperator{\Res}{Res}

\def\v{{\mathbf v}}

\newcommand{\jboxi}[1]{\fbox{\parbox{3mm} { $ \tilde {\imath}_{{#1}}$ } }}

\def\mod{{\rm mod\,}}

\def\+{\mathop{\oplus}}
\def\*{\mathop{\otimes}}

\begin{document}
\title{An observation on highest weight crystals }
\author{ Monica Vazirani}
\maketitle
\section{Introduction }
\label{sec-intro}

As shown in the paper of Stembridge \cite{Stembridge.crystal},
crystal graphs can be characterized  by their local behavior.
In this paper, we observe that a certain local property on
crystals forces a more global property.
In type $A$, this statement says that if a node has  a single parent
and single grandparent, then there is a unique \walk~
from the highest weight node to it.
In other classical types, there is a similar (but necessarily more technical)
statement.
This \walk~ 
is obtained from the associated level $1$ perfect crystal,
$B^{1,1}$. (It is unique unless the Dynkin diagram contains
that of $D_4$ as a subdiagram.)
%

This crystal observation was motivated by representation-theoretic
behavior of the affine Hecke algebra of type $A$, which is
known to be captured by highest weight crystals of type $A^{(1)}$
by the results in 
\cite{Groj.slp}.
As discussed below, the proofs in either setting are straightforward,
and so Grojnowski's theorem linking the two phenomena is not needed.
However, the result is presented here for crystals as one can say
something in all types (Grojnowski's theorem is only in type $A$),
and because  the statement
seems more surprising in the language of crystals than it does for affine Hecke
algebra modules.


\section{Crystals }
\label{sec-crystals}

We begin by reviewing some of the definitions and notation for
crystal graphs, but assume the reader is familiar with crystals and
with root systems. 
For a more comprehensive and complete discussion,
see \cite{Kashiwara.on.crystal.bases}.

\omitt{ wt, epsilon, phi, \e, normal/h.wt. \Lambda_i, G(b)? }
In the following, we fix a root system of finite or affine type.
 $\I$ indexes the simple roots (and the nodes
of the corresponding Dynkin diagram); $P$ is the weight lattice;
$P^*$ is the coroot lattice with canonical pairing $\langle \; , \; \rangle$.
The simple roots are 
$\alpha_i \in P$, and  simple coroots  are $h_i \in P^*$.
The fundamental weights are denoted $\Lambda_i$ and satisfy
$\langle h_i, \Lambda_j \rangle = \delta_{ij}$.
The matrix $[a_{ij}]$ where 
$ a_{ij} = \langle h_i, \alpha_j \rangle$ is the corresponding Cartan
matrix.

A {\it crystal\/} 
is a set of nodes $B$, endowed with the following maps
\begin{align*}
&
\wt: B \to P
\\
&
\epsi: B \to \Z \cup \{ - \infty \}
\\
&
\varphi_i: B \to \Z \sqcup \{ - \infty \}
\\
&
\ei: B \to B \sqcup \{  \zero\}
\\
&
\ffi: B \to B \sqcup \{  \zero\}.
\end{align*}

The maps satisfy the following axioms:
\begin{align*}
&
\varphi_i(\bb) = \epsi(\bb) +   \langle h_i, \wt(\bb) \rangle
\quad \forall i \in \I, \bb \in B.
\\
&
\text{If $\ei \bb \neq \zero$, then } \quad
 \epsi(\ei \bb) = \epsi(\bb) -1,  \quad
	\varphi_i(\ei \bb) = \varphi_i(\bb) +1, \\
& \qquad \qquad \qquad \qquad 
\wt(\ei \bb) = \wt(\bb) + \alpha_i,
\\
&
\text{If $\ffi \bb \neq \zero$, then } \quad
	\epsi(\ffi \bb) = \epsi(\bb) +1,  \quad 
	\varphi_i(\ffi \bb) = \varphi_i(\bb) -1, 
\\
& \qquad \qquad \qquad \qquad 
	\wt(\ffi \bb) = \wt(\bb) - \alpha_i,
\\
&
\text{For $\aA, \bb \in B$, $\aA = \ffi \bb$  if and only if $\bb = \ei \aA$.}
\\
&
\text{ If $\varphi_i(\bb) = - \infty$, then $\ei \bb = \ffi \bb = \zero$.} 
\end{align*}

%
%

Given the crystal data, we can draw the associated {\it crystal graph}.
It is a directed graph with nodes $B$, and $\I$-colored arrows given
by 
$$ \bb \xrightarrow{i} \aA$$
when 
$\aA = \ffi \bb$, or equivalently when  $\bb = \ei \aA$.

In all of the following, we will make the extra assumption that our
crystal $B$ is a highest weight crystal.
\omitt{, [so that in particular it is {\it normal}. ] }
Consequently, we can read the data of
\begin{gather*}
\epsi(\bb) = \max \{n \ge 0 \mid \ei^n \bb \neq \zero \} \\
\varphi_i(\bb) = \max \{n \ge 0 \mid \ffi^n \bb \neq \zero \}, 
\end{gather*}
off of the crystal graph, encoded in the following picture
\begin{gather*}
\bullet
\underbrace{
\xrightarrow{i} \bullet \xrightarrow{i} \cdots 
	\bullet \xrightarrow{i}}_{\epsi(\bb)}
\bb 
\underbrace{ \xrightarrow{i} \bullet \xrightarrow{i} \cdots \xrightarrow{i}
}_{\varphi_i(\bb)}
\bullet .
\end{gather*}
We will also use the notation 
$\eps(\bb) = \sum_{i \in \I} \epsi(\bb) \Lambda_i$.
Thus $\eps(\bb)$ desribes the ``in"-arrows leading to the node $\bb$.

\omitt{
As $\wt(\ei \bb) = \wt(\bb) + \alpha_i$ when $\ei \bb \neq \zero$, we have
\begin{gather}
\label{eq-weight}
\eps(\ei \bb) = \eps(\bb) + \sum_{j \in \I} m_j  \Lambda_j,
\, \text{ where } m_i = -1, 0 \le m_j \le - a_{ij}.
\end{gather}
In general, we have no control over the value $m_j$ takes in the
range $0 \le m_j \le - a_{ij}$ for $j \neq i$.
Below, we will be interested in describing certain cases where
we can force a single $m_j = 1$ and the rest zero. 
}
Below, we will be interested in describing certain cases where
$\eps(\bb) = \Lambda_i$ and $\eps(\ei \bb ) = \Lambda_j$.
(However we will not put any restrictions on ``out"-arrows.)

\subsection{Extra terminology}
\label{sec-terminology}
We introduce some terminology  below.

Let's say that a node $\aA$  is {\it \singular \/}  if
\begin{gather}
\label{eq-star}
\sum_{i \in \I} m_i \le 1, \text{ where } \eps(\aA) = 
\sum_{i \in \I} m_i  \Lambda_i.
\end{gather}
Notice equation \eqref{eq-star} implies there is at most one
$i \in \I$ such that $\ei \aA \neq \zero$.  In particular, highest
weight nodes satisfy \eqref{eq-star}.
In the crystal graph, we picture \singular~ nodes as having a single
``in"-arrow leading to it (and any arrow preceeding that one carries a 
different color), but there is no restriction on its ``out"-arrows.

If $\ei (\aA) = \bb$ for some $i$, we shall say $\bb$ is a
{\it parent\/} of $\aA$.
We will define {\it ancestor\/} inductively by saying parents are ancestors
and parents of ancestors are also ancestors.

\section{Kashiwara's Theorem for Highest Weight Crystals }
\label{sec-kashiwara}

In all the following theorems, we fix a 
root system and
assume $B$ is a fixed highest
weight crystal of that type.

The crystal graph $B$ comes from an integrable highest weight module $V$ of
the associated Lie algebra or quantum enveloping algebra.
We appeal to theorems of Kashiwara that ensure the existence
of a global basis 
$\{ G(\bb) \mid \bb \in B \}$ of $V$. 
In the following $\eei$ will denote a 
Chevalley generator, and
$\eei^{(m)}$ its divided power.

We first give a remark (in Section 5) of \cite{Kashiwara.global}
as the following useful lemma.
One should compare it to the statement 
 $\wt(\ei \bb) = \wt(\bb) + \alpha_i$. 
\begin{lemma}
\label{lemma-eps}
When $\ei \bb \neq \zero$,
\begin{gather}
\label{eq-weight}
\eps(\ei \bb) = \eps(\bb) + \sum_{j \in \I} m_j  \Lambda_j,
\, \text{ where } m_i = -1,\quad 0 \le m_j \le - a_{ij}.
\end{gather}
\end{lemma}
In general, we have no control over the value $m_j$ takes in the
range $0 \le m_j \le - a_{ij}$ for $j \neq i$.
Below, we will be interested in describing certain cases where
we can force a single $m_j = 1$ and the rest zero. 
In other words, we want that $\eps(\bb) = \Lambda_i$
and $\eps(\ei \bb ) = \Lambda_j$.

We list some immediate corollaries to this lemma.
\begin{corollary}
\label{cor-zero}
Let $ \bb \in B$ and suppose $a_{ij} =  0$. 
Then $\ej \bb = \zero \implies \ej \ei \bb  = \zero$.
\end{corollary}
\begin{corollary}
\label{cor-parent}
Let $\aA, \bb \in B$ both be \singular~ nodes, and suppose 
$\bb$ is  a parent of $\aA$, with $\bb  = \ei \aA$.
Then $\ej \bb \neq \zero \implies a_{ij} <  0$. 
\end{corollary}


\begin{theorem}[\cite{Kashiwara.global}]
\label{thm-kash}
Let $\bb \in B$ and suppose $\ei^m \bb \neq \zero$ but 
 $\ei^{m+1} \bb = \zero$.
Then
$$\eei^{(m)} G( \bb) = G( \ei^m \bb) \quad \text{ and }
\quad \eei^{(m+1)} G( \bb) = 0.$$
\end{theorem}

As a corollary to this theorem, employing the Serre relations,
we can deduce several properties of \singular~ nodes.
We review the Serre relations below.

Fix $i, j \in \I, i \neq j$.  Let $\ell = 1 - \langle h_i, \alpha_j \rangle
= 1 - a_{ij}$. 
Then
\begin{gather}
\label{eq-serre}
\sum_{k=0}^\ell \eei^{(k)} \eej \eei^{(\ell -k)} = 0.
\end{gather}
\begin{corollary}
\label{cor-serre}
\begin{enumerate}
\item
\label{case-0}
Suppose $a_{ij} = 0$. Then $\ej \bb = \zero, \ei^2 \bb = \zero
	 \implies \ej(\ei \bb) = \zero$.
\item
\label{case-1}
Suppose $a_{ij} = -1$. Then $\ej \bb = \zero, \ei^2 \bb = \zero
 \implies \ei(\ej \ei \bb) = \zero$.
\item
\label{case-2}
Suppose $a_{ij} = -2$.
Then $\ej \bb = \zero, \ei^2 \bb = \zero, \ej^2 \ei \bb = \zero
 \implies \ei(\ei \ej \ei \bb) = \zero$.
Also $\ei b = \zero, \ej^2 \bb = \zero \implies \ei^3 \ej \bb = \zero$.

If in addition  $a_{ji} = -1$, then
 $\ej b = \zero, \ei^2 \bb = \zero \implies
 \ej^2 \ei \bb = \zero $ and $\ej (\ei \ej \ei \bb) = \zero$.
\end{enumerate}
\end{corollary}
\begin{proof}
\eqref{case-0} 
This follows directly from 
Corollary \ref{cor-zero}, which is a stronger statement.
(We note one may also prove this using Theorem \ref{thm-kash} 
in a manner similar to the subsequent cases.)

\eqref{case-1} 
From the Serre relations for $a_{ij} = -1$, we know that 
$(\eei^{(2)} \eej - \eei \eej \eei + \eej \eei^{(2)})(G(\bb)) = 0$.
Applying Theorem \ref{thm-kash}, $\ej \bb = \zero \implies \eej G(\bb) = 0$
and 
 $\ei^2 \bb = \zero$ implies both that $ \eei^{(2)} G(\bb) = 0$
and $\eei G(\bb) = G(\ei \bb)$.
Hence we get $0= \eei \eej \eei G(\bb) = \eei \eej G(\ei \bb) $.

Kashiwara's equation 5.3.8 in \cite{Kashiwara.global} 
gives $\eei G(\bb)$ as a linear combination of $G(\ei \bb)$ and
$G(\bb')$ where $\varphi_k(\bb') \le \varphi_k(\bb)$ for all $k \in I$.
Iterating this, we get that $0 = \eei \eej G(\ei \bb)$ is a linear
combination of $ G(\ei \ej \ei \bb)$ and terms $G(\bb')$. 
It is straightforward (using equation \eqref{eq-weight})
 to show the restrictions on $\bb'$ can only
be satisfied if $\epsi(\ej \ei \bb) \le -a_{ij} -1 = 0$.
But this forces $\ei \ej \ei \bb = \zero$.  In the case there are
no such $\bb'$, we then get $G(\ei \ej \ei \bb) = 0$, so again
$\ei \ej \ei \bb = \zero$.

\eqref{case-2} 
The conditions on $\bb$ 
give us $\eej G(\bb) = 0$, $\eei G(\bb) = G(\ei \bb)$, $\eei G(\ei \bb) = 0$,
and $\eej \eei G(\bb) = G(\ej \ei \bb).$
The Serre relations imply 
$0 = \eei^{(2)} \eej \eei G(\bb) = \eei^{(2)} G(\ej \ei \bb) $ 
which implies $\ei^2 \ej\ei \bb = \zero$.
(In particular, this also implies
$\eei \eej \eei G( \bb) = G( \ei \ej \ei \bb)$.)
For the second case, we get $0 = \eei^{(3)} \eej G(\bb) =\eei^{(3)} G(\ej \bb)$
so that $\ei^3 \ej \bb = \zero$.

For the final statement, the proof of the first implication follows
immediately from equation \eqref{eq-weight}.  For the second, as $a_{ji}
= -1$,
$\eej \eei \eej \eei G(\bb) = \eej \eei \eej G(\ei \bb) = 
 \eej^{(2)} \eei  G(\ei \bb) +  \eei \eej^{(2)}  G(\ei \bb) =  0$.
So $ 0 = \eej \eei \eej G(\ei \bb) = \eej G( \ei \ej \ei \bb)$,
yielding $\ej \ei \ej \ei \bb = \zero$.

\omitt{
For the other cases, the proofs are similar: the hypotheses ensure
all but one of the terms in the Serre relation (applied to $G(\bb)$) vanish,
and applying Theorem \ref{thm-kash} to the remaining term gives us vanishing
in the crystal. 
}

\end{proof}
We remark that there are similar statements for $a_{ij} = -3, -4$, but they
do not translate into interesting statements about \singular~ nodes as
the other cases do in Theorem \ref{thm-type} below.

Case \eqref{case-0} of  Corollary \ref{cor-serre} says that if $a_{ij} = 0$
and we
see
 $$\bullet \xrightarrow{i} \bb,\text{ we do not see } \bullet \xrightarrow{j}
\bullet \xrightarrow{i} \bb.$$
Compare 
this  with the fact that $a_{ij} = 0$ means 
that in the Dynkin diagram we see
\newline
\begin{center}
\begin{picture}(26,20)(-5,-10)
\multiput( 0,0)(20,0){2}{\circle{6}}
\put(0,-5){\makebox(0,0)[t]{$j$}}
\put(20,-5){\makebox(0,0)[t]{$i$}}
\end{picture}
$\,$ and not $\,$  
\begin{picture}(26,20)(-5,-10)
\multiput( 0,0)(20,0){2}{\circle{6}}
\multiput( 3,0)(20,0){1}{\line(1,0){14}}
\put(0,-5){\makebox(0,0)[t]{$j$}}
\put(20,-5){\makebox(0,0)[t]{$i$}}
\end{picture} 
\end{center}
in the crystal.
Similarly, when 
$a_{ij} = -1$
and we see
 $$\bullet \xrightarrow{j} \bullet \xrightarrow{i} \bb,
\text{ we do not see }
 \bullet \xrightarrow{i} \bullet \xrightarrow{j} \bullet \xrightarrow{i} \bb.$$
Compare this to the fact that when $a_{ij} = -1$ 
we see
\begin{picture}(26,26)(-5,-10)
\multiput( 0,0)(20,0){2}{\circle{6}}
\multiput( 3,0)(20,0){1}{\line(1,0){14}}
\put(0,-5){\makebox(0,0)[t]{$j$}}
\put(20,-5){\makebox(0,0)[t]{$i$}}
\end{picture}
in the Dynkin diagram
but
\newline
\begin{center}
not $\,$
\begin{picture}(26,20)(-5,-10)
\multiput( 0,0)(20,0){2}{\circle{6}}
\multiput(2.85,-1)(0,2){2}{\line(1,0){14.3}}
\put(0,-5){\makebox(0,0)[t]{$i$}}
\put(20,-5){\makebox(0,0)[t]{$j$}}
\put( 6, 0){\makebox(0,0){$<$}}
\put(14, 0){\makebox(0,0){$>$}}
\end{picture}
 $\,$ nor $\,$
\begin{picture}(26,20)(-5,-10)
\multiput( 0,0)(20,0){2}{\circle{6}}
\multiput( 2.85,-1)(0,2){2}{\line(1,0){14.3}} 
\put(10,0){\makebox(0,0){$>$}}
\put(0,-5){\makebox(0,0)[t]{$j$}}
\put(20,-5){\makebox(0,0)[t]{$i$}}
%
\end{picture}
 $\,$ nor $\,$
\begin{picture}(26,20)(95,-10)
\multiput(100,0)(20,0){2}{\circle{6}}
\multiput(102.85,-1)(0,2){2}{\line(1,0){14.3}} 
\put(110,0){\makebox(0,0){$<$}}
\put(100,-5){\makebox(0,0)[t]{$i$}}
\put(120,-5){\makebox(0,0)[t]{$j$}}
%
\end{picture}
\end{center}
which we should associate to
\begin{picture}(26,20)(-5,-10)
\multiput( 0,0)(20,0){2}{\circle{6}}
\put(0,-5){\makebox(0,0)[t]{$i$}}
\put(20,-5){\makebox(0,0)[t]{$j$}}
\put(10,3){\oval(20,7)[t]}
\put(10,-3){\oval(20,7)[b]}
\end{picture}
in the former case, and to the folding of
\begin{picture}(26,40)(-5,-10)
\multiput( 0,0)(20,0){2}{\circle{6}}
\multiput(20,20)(60,0){1}{\circle{6}}
\multiput( 3,0)(20,0){1}{\line(1,0){14}}
\multiput(20,3)(60,0){1}{\line(0,1){14}}
\put(0,-5){\makebox(0,0)[t]{$i$}}
\put(20,-5){\makebox(0,0)[t]{$j$}}
\put(25,20){\makebox(0,0)[l]{$i$}}
\end{picture}
%
in the latter cases.
In Theorem \ref{thm-type} below, we shall see that requiring certain
\singular ity conditions on nodes forces the colors on their in-arrows
to 
behave as  
a directed ``path" or \walk~  would on the Dynkin diagram, as suggested above.
%
Choosing $\aA, \bb$ with $\eps(\aA) = \Lambda_i$, $\eps(\bb) = \Lambda_j$
and $\bb$ the parent of $\aA$ puts an ``orientation" on the Dynkin
diagram. As the Dynkin diagram's  vertices correspond to arrows
in the crystal, we really are making a statement about a graph
{\it dual\/} to the Dynkin diagram.
It turns out the correct notion of duality in this setting is
exactly captured in an associated level $1$ perfect crystal.


Below, we recap, case by case, the consequences of
Corollary \ref{cor-serre}
on all of the ancestors of a singular node $\aA$ whose parent
is also \singular.
We  describe all \walks~ on the crystal,
from the highest weight node $\v$ to $\aA$.
These \walks~ are described exactly by \walks~ on the level $1$ perfect
crystal $B^{1,1}$.
These crystals are displayed in the body of the proof as well
as in the appendix.
(In type $A_n^{(1)}$ we also need the
perfect crystal $B^{n,1}$ obtained by reversing all arrows in $B^{1,1}$.
In type $A_1^{(1)}$ we require the grandparent to be \singular~ as well.)

A necessary, but not sufficient, condition for both a node and its
parent to be \singular~ is that it has the form
$  \ff{i_1} \ff{i_2} \cdots \ff{i_k} \v,$
where $\xrightarrow{i_1} \xrightarrow{i_2}\cdots \xrightarrow{i_k}$
is a consecutive sequence   of arrows in  $B^{1,1}$.
(The theorem also describes which nodes of this form are not \singular.)
This means that we can give a case by case description of the node's
ancestors, but a global statement about the \walks~ from $\v$ to $\aA$.
The local nature of \singular ity means that the result in affine
type follows from that in finite type (in small rank), and so we structure the
statements and proofs of the following theorem accordingly.

\begin{theorem}
\label{thm-global}
Let $B$ be a highest weight crystal with highest weight node $\v$
 of type
$A_n, n \ge 1$,
$A_n^{(1)}, n \ge 2$,
$A_{2n}^{(2)}, n \ge 2$,
$A_{2n}^{(2) \dagger}, n \ge 2$,
$A_{2n-1}^{(2)}, n \ge 3$,
$B_{n}, n \ge 2$,
$B_{n}^{(1)}, n \ge 3$,
 $C_n, n \ge 2$,
 $C_n^{(1)}, n \ge 2$,
$D_{n}, n \ge 4$,
$D_{n}^{(1)}, n \ge 4$,
$D_{n+1}^{(2)}, n \ge 2$.
Suppose $\aA \in B$ is a \singular~ node with  \singular~ parent.
Then 
 $$ \aA =  \ff{i_1} \ff{i_2} \cdots \ff{i_k} \v,$$
only when
$\xrightarrow{i_1} \xrightarrow{i_2}\cdots \xrightarrow{i_k}$
is a consecutive sequence   of arrows in  the level $1$
perfect crystal  $B^{1,1}$ (or $B^{n,1}$ in type $A$)
 of appropriate type,
omitting $0$-arrows in finite type.
If the Dynkin diagram does not contain that of $D_4$ as a subdiagram,
then this sequence is unique.
(In type 
$A_1^{(1)}$, we get the same conclusion if we also require
$\aA$ also have \singular~ grandparent.)
\end{theorem}
\begin{theorem}
\label{thm-type}
\begin{enumerate}
\item
\label{case-A}
Let $B$ be a highest weight crystal of type $A_n, n \ge 1$
or of type $A_n^{(1)}, n \ge 2$.

Suppose $\aA, \bb \in B$ are both \singular~ nodes with $\bb  = \ei \aA$.
Then all ancestors of $\aA$ are \singular.
There is a {\it unique\/} \walk~ (on the directed graph) from
the highest weight node  $\v \in B$ to $\aA$, given by 
$\aA = \ffi \ff{i \pm 1} \ff{i \pm 2} \cdots \ff{i \pm k} \v,$
where subscripts are taken $\mod n$.
\item
\label{case-Aone}
Let $B$  be  a highest weight crystal of type 
type $A_1^{(1)}$.

Suppose $\aA, \bb, \cc \in B$ are all \singular~ nodes with $\bb  = \ei \aA$,
$\cc  = \ej \bb$.  (Necessarily, $i \neq j$.)
Then all ancestors of $\aA$ are \singular.
There is a {\it unique\/}  \walk~ from
the highest weight node  $\v $ 
to $\aA$, given by 
$\aA = \ffi \fj \ffi  \fj \cdots \v.$
\item
\label{case-C}
Let $B$ be a highest weight crystal of type 
$C_n, n \ge 2$.

Suppose $\aA, \bb \in B$ are both \singular~ nodes with $\bb  = \ei \aA$.
Then all ancestors of $\aA$ are \singular.
There is a {\it unique\/}  \walk~
 from the highest weight node  $\v$ 
to $\aA$, given by 
the following possibilities:
\omitt{ (which are also describable as
 a sequence of consecutive arrows on the level $1$ perfect crystal
$B^{1,1}$ of affine type $C_n^{(1)}$, omitting the $0$-arrow):
}
\begin{enumerate}
\item
$\aA = \ffi \ff{i \pm 1}  \cdots \v.$
\item
$\aA = \ffi \ff{i + 1} \ff{i + 2} \cdots \ff{n -1} \ff{n} \ff{n-1}
\ff{n-2} \cdots \v.$
\end{enumerate}
\item
\label{case-Bfin}
Let $B$ be a highest weight crystal of type $B_n, n \ge 2$.

Suppose $\aA, \bb \in B$ are both \singular~ nodes with $\bb  = \ei \aA$.
Then all but one of the ancestors of $\aA$ are \singular.
If there is a non-\singular~ ancestor $\cc$, it satisfies $\eps(\cc) = 
2 \Lambda_n$.

There is a {\it unique\/}  \walk~ 
 from the highest weight node  $\v$ 
to $\aA$, given by 
the following possibilities:
	\begin{enumerate}
\item
$\aA = \ffi \ff{i \pm  1} \ff{i \pm  2}  \cdots \v.$
\item 
$\aA = \ffi \ff{i + 1}  \cdots \ff{n-1} \ff{n} \ff{n} \ff{n-1}
\ff{n-2} \cdots  \v.$
	\end{enumerate}

\item
\label{case-Dfin}
Let $B$ be a highest weight crystal of type $D_n, n \ge 4$.

Suppose $\aA, \bb \in B$ are both \singular~ nodes with $\bb  = \ei \aA$.
Then all but one of the ancestors of $\aA$ are \singular.
If there is a non-\singular~ ancestor $\cc$, it satisfies
$\eps(\cc) = \Lambda_{n-1} +  \Lambda_n$.

There are at most two  \walks~
 from the highest weight node  $\v$ 
to $\aA$, given by 
the following possibilities.
\omitt{
(which are also describable as
 a sequence of consecutive arrows on the level $1$ perfect crystal
$B^{1,1}$ of affine type $D_n^{(1)}$, omitting the $0$-arrow).
}

Below we use the notation
$\ff{n-1, n}$ to stand for  either
$ \ff{n-1} \ff{n} $ or $ \ff{n} \ff{n-1} $.
  (Of course, it is possible
the path truncates before giving both terms of $\ff{n-1, n}$ 
which would give a unique  \walk.)
	\begin{enumerate}
\item
$\aA = \ffi \ff{i \pm  1}  \cdots \v$
\item 
$\aA = \ffi \ff{i + 1}  \cdots \ff{n-2}  \ff{n-1, n}  \ff{n-2}
\cdots  \v$  (In this case, we get two  \walks.)
	\end{enumerate}
%
\item
\label{case-therest.unique}
Let $B$ be a highest weight crystal of type
 $C_n^{(1)}, n \ge 2$.
$A_{2n}^{(2)}, n \ge 2$,
$A_{2n}^{(2) \dagger}, n \ge 2$,
or
$D_{n+1}^{(2)}, n \ge 2$.

Suppose $\aA, \bb \in B$ are both \singular~ nodes with $\bb  = \ei \aA$.
Then ancestors $\cc$ of $\aA$ are either  \singular~
or they satisfy
 $\eps(\cc) =2\Lambda_{n}$ in types
		$A_{2n}^{(2) \dagger}$, $D_{n+1}^{(2)}$; 
 $\eps(\cc) =2\Lambda_{0}$ in types $A_{2n}^{(2)}$, $D_{n+1}^{(2)}$. 

In all cases, there is a {\it unique\/}  \walk~
 from the highest weight node  $\v$ 
to $\aA$, given by 
the following possibilities:
 $$ \aA =  \ff{i} \ff{i_2} \cdots \ff{i_k} \v,$$
where $\xrightarrow{i} \xrightarrow{i_2}\cdots \xrightarrow{i_k}$
is a consecutive sequence   of arrows in  the level $1$
perfect crystal  $B^{1,1}$ of appropriate type.

\item
\label{case-therest.infinite}
Let $B$ be a highest weight crystal of type
$D_{n}^{(1)}, n \ge 4$,
$A_{2n-1}^{(2)}, n \ge 3$,
or 
$B_{n}^{(1)}, n \ge 3$.

Suppose $\aA, \bb \in B$ are both \singular~ nodes with $\bb  = \ei \aA$.
Then ancestors $\cc$ of $\aA$ are either  \singular~
or they satisfy
 $\eps(\cc) =2\Lambda_{n}$ in type
	$B_{n}^{(1)}$;
$\eps(\cc) = \Lambda_{n-1} +  \Lambda_n$, in type  
	$  D_{n}^{(1)}$;
$\eps(\cc) = \Lambda_{1} +  \Lambda_0$ in types  $ D_{n}^{(1)}, B_{n}^{(1)}$.

Walks
 from the highest weight node  $\v \in B$ to $\aA$, described by 
the following (infinite)  possibilities.
 $$ \aA =  \ff{i} \ff{i_2} \cdots \ff{i_k} \v,$$
where $\xrightarrow{i} \xrightarrow{i_2}\cdots \xrightarrow{i_k}$
is a consecutive sequence   of arrows in  $B^{1,1}$.
\end{enumerate}
\end{theorem}
%
We remark that in cases not included above, such as exceptional types,
or type
$A_{2}^{(2)}$,  that requiring a certain number of consecutive \singular~
nodes either gives many possible complicated
\walks~ from  the highest weight node or none at all. 
At the end of this paper we have a short discussion regarding type
$E_6$.

\omitt{
We also note that each of the paths back to the highest weight node
correspond to following a path on the corresponding Dynkin diagram,
placing a ``u-turn" at each inward facing arrow, and a mirror after each
trivalent branching. [wording?]
}

\begin{proof}
\eqref{case-A} [$A_n$, $A_n^{(1)}$]
We have $\eps(\aA) = \Lambda_i$, and either $\bb = \v$ or $\eps(\bb) 
= \Lambda_j$ with $j$ connected to $i$ in the Dynkin diagram by Corollary
\ref{cor-parent}.  In this case $j = i \pm 1$, taking 
$j \mod n$ if
necessary.  Applying this corollary again, $\ee{k} (\ej \bb) = \zero$
unless $k = j \pm 1$.  By case \eqref{case-1} of Corollary \ref{cor-serre},
$\zero = \ei(\ej \ei \aA) = \ei \ej \bb$, so we must have $k = i \pm 2$
and either $ \ej \bb = \v$ or $\eps( \ej \bb) = \Lambda_k$.  Hence
we can inductively apply this argument to the pair $\bb$ and $\ej \bb$.  
As $B$ is a highest weight crystal, this process must eventually
terminate at $\ee{i \pm m} \cdots \ee{i \pm 1} \ei \aA = \v$ which is equivalent
to 
$\aA = \ffi \ff{i \pm 1} \ff{i \pm 2} \cdots \ff{i \pm m} \v.$

The above sequence of consecutively colored arrows exactly corresponds
to a sequence of arrows on the following perfect crystals.
$$
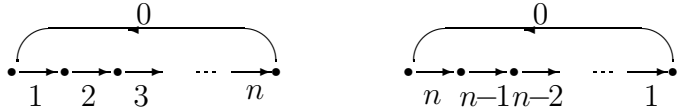
\begin{picture}(150,48)(-10,-12)
\multiput( 0,0)(20,0){3}{\circle*{3}}
\multiput( 3,0)(20,0){3}{\vector(1,0){14}}
\multiput(70,0)(3,0){3}{\line(1,0){1}} 
\put(100,0){\circle*{3}}
\put(83,0){\vector(1,0){15}}
\put(51,4){\oval(98,25)[t]}
\put(47,16){\vector(-1,0){3}}
\put(9,-5){\makebox(0,0)[t]{$1$}}
\put(29,-5){\makebox(0,0)[t]{$2$}}
\put(49,-5){\makebox(0,0)[t]{$3$}}
\put(92,-5){\makebox(0,0)[t]{$n$}}
\put(50,26){\makebox(0,0)[t]{$0$}}
\end{picture}
\begin{picture}(150,48)(-10,-12)
\multiput( 0,0)(20,0){3}{\circle*{3}}
\multiput( 3,0)(20,0){3}{\vector(1,0){14}}
\multiput(70,0)(3,0){3}{\line(1,0){1}} 
\put(100,0){\circle*{3}}
\put(83,0){\vector(1,0){15}}
\put(51,4){\oval(98,25)[t]}
\put(47,16){\vector(-1,0){3}}
\put(9,-7){\makebox(0,0)[t]{$n$}}
\put(29,-5){\makebox(0,0)[t]{$n\! \!-\! \!1$}}
\put(49,-5){\makebox(0,0)[t]{$n\! \!-\! \!2$}}
\put(92,-5){\makebox(0,0)[t]{$1$}}
\put(50,26){\makebox(0,0)[t]{$0$}}
\end{picture}
$$

As we only care about the arrow labels, we omit the node labels
that are usually also pictured in the crystals.

The reader should compare the above perfect crystals
to  the Dynkin diagrams
$$
A_n:
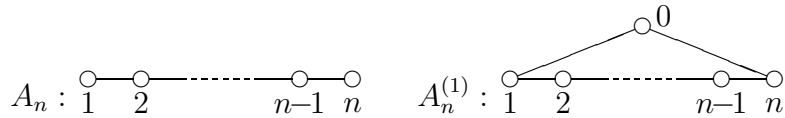
\begin{picture}(106,35)(-5,-10)
\multiput( 0,0)(20,0){2}{\circle{6}}
\multiput(80,0)(20,0){2}{\circle{6}}
\multiput( 3,0)(20,0){2}{\line(1,0){14}}
\multiput(63,0)(20,0){2}{\line(1,0){14}}
\multiput(39,0)(4,0){6}{\line(1,0){2}}
\put(0,-5){\makebox(0,0)[t]{$1$}}
\put(20,-5){\makebox(0,0)[t]{$2$}}
\put(80,-5){\makebox(0,0)[t]{$n\!\! -\!\! 1$}}
\put(100,-7){\makebox(0,0)[t]{$n$}}
\end{picture}
\qquad 
A_n^{(1)}:
\begin{picture}(106,35)(-5,-10)
\multiput( 0,0)(20,0){2}{\circle{6}}
\multiput(80,0)(20,0){2}{\circle{6}}
\put(50,20){\circle{6}}
\multiput( 3,0)(20,0){2}{\line(1,0){14}}
\multiput(63,0)(20,0){2}{\line(1,0){14}}
\multiput(39,0)(4,0){6}{\line(1,0){2}}
\put(2.78543,1.1142){\line(5,2){44.429}}
\put(52.78543,18.8858){\line(5,-2){44.429}}
\put(0,-5){\makebox(0,0)[t]{$1$}}
\put(20,-5){\makebox(0,0)[t]{$2$}}
\put(80,-5){\makebox(0,0)[t]{$n\!\! -\!\! 1$}}
\put(100,-7){\makebox(0,0)[t]{$n$}}
\put(55,20){\makebox(0,0)[lb]{$0$}}
\end{picture}
$$
and  consider the discussion below Corollary \ref{cor-serre}.  Observe that 
arrows being consecutive in the perfect crystal correspond to
vertices being adjacent in the Dynkin diagram.

\eqref{case-Aone}[$A_1^{(1)}$]
We note that only in this case do  we require {\it three\/} consecutive
\singular~ nodes.
As above, we necessarily have $\eps(\aA) = \Lambda_i$, $\eps(\bb) = 
\Lambda_j$, $\eps(\cc) = \Lambda_i$ (or $\cc = \v$)
 by Corollary \ref{cor-parent}.
If we can show $\eps(\ei \cc) = \Lambda_j$ or that $\ej \ei \cc = \zero$ 
(forcing $\ei \cc = \v$), we will be done by a similar induction as used above.
We already are given $\ei^2 \cc = \zero$.  Observe $\ej^2 \ei \cc =
\ej^2 \ei  \ej \bb = \zero$ by case \ref{case-2} of Corollary \ref{cor-serre}.
Again, the reader can compare this statement to tracing a path on
the Dynkin diagram
$$A_1^{(1)}: 
\begin{picture}(26,25)(-5,-10)
\multiput( 0,0)(20,0){2}{\circle{6}}
\multiput(2.85,-1)(0,2){2}{\line(1,0){14.3}}
\put(0,-5){\makebox(0,0)[t]{$0$}}
\put(20,-5){\makebox(0,0)[t]{$1$}}
\put( 6, 0){\makebox(0,0){$<$}}
\put(14, 0){\makebox(0,0){$>$}}
\end{picture}
\quad\text{which  is again captured in the perfect crystal} \quad
\begin{picture}(26,25)(-5,-10)
\multiput( 0,0)(20,0){2}{\circle{6}}
\put(8,-10){\makebox(0,0)[t]{$1$}}
\put(12,15){\makebox(0,0)[t]{$0$}}
\put(9,-7){\vector(1,0){5}}
\put(11,6){\vector(-1,0){5}}
\put(10,3){\oval(20,7)[t]}
\put(10,-3){\oval(20,7)[b]}
\end{picture} \, \; .
$$

\eqref{case-C}[$C_n$]
This proof is 
similar to that of
case \eqref{case-A}.  
We need only consider
 the case that $\aA$ and $\bb = \ee{n-1} \aA$ are \singular~ with 
$\eps(\aA) = \Lambda_{n-1}, \eps(\bb) = \Lambda_n$.
Let $\cc = \ee{n} \bb$.  We claim either $\cc  = \v$ or $\eps(\cc) = 
\Lambda_{n-1}$.  By Corollary \ref{cor-parent}, we know
$\ee{k} \cc = \zero$ unless $k = n-1$.   Because $a_{n-1,n} = -2,$
by case \ref{case-2} of Corollary \ref{cor-serre}, we know
$\ee{n-1}^2 \ee{n} \ee{n-1} \aA = \zero$.  This gives the claim.
Now the induction proceeds just as in type $A$.

Again, we draw the Dynkin diagram
$$
 C_n:
\begin{picture}(126,20)(10,-10)
\multiput( 20,0)(20,0){2}{\circle{6}}
\multiput(100,0)(20,0){2}{\circle{6}}
\multiput(23,0)(20,0){2}{\line(1,0){14}}
\put(83,0){\line(1,0){14}}
\multiput(102.85,-1)(0,2){2}{\line(1,0){14.3}} 
\multiput(59,0)(4,0){6}{\line(1,0){2}} 
\put(110,0){\makebox(0,0){$<$}}
\put(20,-5){\makebox(0,0)[t]{$1$}}
\put(40,-5){\makebox(0,0)[t]{$2$}}
\put(100,-5){\makebox(0,0)[t]{$n\!\! -\!\! 1$}}
\put(120,-7){\makebox(0,0)[t]{$n$}}
\end{picture}
$$
and show the perfect crystal of type $C_n^{(1)}$ with the $0$-arrow 
removed, which is suggestive  of picturing the double arrow as a folding.
(Note that we recover the same graph reversing orientation of all arrows.)
$$
\qquad 
\begin{picture}(126,35)(10,-10)
\multiput( 20,0)(20,0){2}{\circle{6}}
\put(120,8){\circle{6}}
\put(103,0){\line(5,3){14}}
\put(103,15){\line(5,-3){14}}
\put(104,0){\vector(4,3){8}}
\put(111,10){\vector(-4,3){5}}
\put(83,0){\vector(1,0){8}}
\multiput(35,15)(20,0){2}{\vector(-1,0){8}}
\put(95,15){\vector(-1,0){8}}
\multiput(23,0)(20,0){2}{\vector(1,0){8}}
\multiput(100,0)(20,0){1}{\circle{6}}
\multiput(23,0)(20,0){2}{\line(1,0){14}}
\put(83,0){\line(1,0){14}}
\multiput( 20,15)(20,0){2}{\circle{6}}
\multiput(100,15)(20,0){1}{\circle{6}}
\multiput(23,15)(20,0){2}{\line(1,0){14}}
\put(83,15){\line(1,0){14}}
\multiput(59,0)(4,0){6}{\line(1,0){2}} 
\multiput(59,15)(4,0){6}{\line(1,0){2}} 
\put(20,-5){\makebox(0,0)[t]{$1$}}
\put(40,-5){\makebox(0,0)[t]{$2$}}
\put(100,-5){\makebox(0,0)[t]{$n\!\! -\!\! 1$}}
\put(120,-1){\makebox(0,0)[t]{$n$}}
\put(20,27){\makebox(0,0)[t]{$1$}}
\put(40,27){\makebox(0,0)[t]{$2$}}
\put(100,27){\makebox(0,0)[t]{$n\!\! -\!\! 1$}}
\end{picture} 
$$
Note that the conclusions (a),(b) can
also be expressed as
$\aA = \ffi \ff{i \pm 1} \ff{i \pm 2} \cdots \ff{i \pm k} \v,$
so long as subscripts are taken $\mod 2n$, and one sets $\ff{n+m} := 
\ff{n-m}$ for $0 < m < n$.

\eqref{case-Bfin}[$B_n$]
We need only consider
 the case that $\aA$ and $\bb = \ee{n-2} \aA$ are \singular~ with 
$\eps(\aA) = \Lambda_{n-2}, \eps(\bb) = \Lambda_{n-1}$.
Otherwise it reduces to case \eqref{case-A}.
Let $\cc = \ee{n-1} \bb$.  We claim either $\cc  = \v$;
$\eps(\cc) = \Lambda_{n}$ in which case $\ee{n} \cc = \v$; or 
$\eps(\cc) = 2 \Lambda_{n}$, in which case $\cc$ is {\it not\/}
\singular, but  both $\ee{n} \cc$ and
$\ee{n}^2 \cc$ are \singular,
and $\ee{n-1} \ee{n}^2 \cc$ is either \singular~ or $\zero$. 

By Corollary \ref{cor-parent}, we know
$\ee{k} \cc = \zero$ unless $k = n$ or $n-2$. But case \eqref{case-A} of
this theorem rules out the latter.   Because $a_{n,n-1} = -2,$
by case \ref{case-2} of Corollary \ref{cor-serre} we know
$\ee{n}^3 \ee{n-1} \bb = \zero$, showing the first part of the claim.
Now suppose $\eps(\cc) =  \Lambda_{n}$.
That means $\ee{n}^2 \cc = \zero$. Further, $\ee{n-1} \ee{n} \cc
= \ee{n-1} \ee{n} \ee{n-1} \bb = \zero$ by case \eqref{case-1} of
Corollary \ref{cor-serre} as $a_{n-1,n} = -1$. 
By Corollary \ref{cor-parent}, $\ee{k} \ee{n} \cc = \zero$ for all
$k \neq n-1$, showing $\ee{n} \cc  = \v$ as the crystal $B$ has a unique
highest weight node.

Next suppose $\eps(\cc) = 2 \Lambda_{n}$.
In particular, notice that $\cc$ is not \singular.
\omitt{
We note that $\ee{n} \cc$ is \singular, because first
$\ee{n-1} \ee{n} \cc = \ee{n-1} \ee{n} \ee{n-1} \bb = \zero$ by
case \eqref{case-1} of Corollary \ref{cor-serre}.
Secondly,
by \eqref{case-0},
$k \neq n-1$,
$\ee{k} \ee{n} \cc =  \zero$.
Hence $\eps( \ee{n} \cc) = \Lambda_n$.
}
For   $k \neq n-1, n$,
\omitt{
as $\ee{k} \cc = \zero$, we know $0 = \eee{n}^{(2)} \eee{k} G(\cc)
=  \eee{k}  \eee{n}^{(2)} G(\cc) =  \eee{k}  G( \ee{n}^2 \cc)$
by Theorem \ref{thm-kash}, showing
$\ee{k} \ee{n}^2 \cc =  \zero$.
}
we know $0 = \epsof{k}(\cc) =  \epsof{k}(\ee{n} \cc) = 
\epsof{k}( \ee{n}^2 \cc)$ by Lemma \ref{lemma-eps}.
As above, we still have $\ee{n-1} \ee{n} \cc = \zero$, so that
$\ee{n} \cc$ is \singular.

To show $\ee{n}^2 \cc$ is \singular, we need only show $\ee{n-1}^2 \ee{n}^2 
\cc = \zero$.
Note $\eee{n-1}^{(2)} G(\ee{n}^2 \cc ) =
\eee{n-1}^{(2)} G(\ee{n}^2 \ee{n-1} \bb)
=  \eee{n-1}^{(2)}\eee{n-1}^{(2)}\eee{n-1} G(\bb)
=$ 
$(\eee{n-1} \eee{n}^{(2)} - \frac{1}{[2]![2]!} \eee{n} \eee{n-1} \eee{n})
\eee{n-1}^{(2)}  G(\bb)
+ 
(\frac{1}{[2]!} \eee{n-1} \eee{n}\eee{n-1}^{(2)} - 
\frac{[3]}{[2]![2]!} \eee{n} \eee{n-1}\eee{n-1}^{(2)}) \eee{n} G(\bb)
= 0$ 
by the Serre relations and  Theorem \ref{thm-kash}.
Hence  $\ee{n-1}^2 \ee{n}^2 \cc = \zero$.


\omitt{
Let $B$ be a highest weight crystal of type $B_n, n \ge 2$, 
$A_{2n}^{(2)}, n \ge 2$,
$D_{n+1}^{(2) }, n \ge 2$.
or of type $A_{2n}^{(2) \dagger}, n \ge 2$.
}
We again show the Dynkin diagram, along with the perfect crystal
with $0$-arrows removed.
$$B_n: \; \,
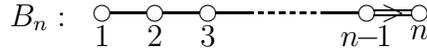
\begin{picture}(126,30)(-5,-5)
\multiput(0,0)(20,0){3}{\circle{6}}
\multiput(100,0)(20,0){2}{\circle{6}}
\multiput( 3,0)(20,0){3}{\line(1,0){14}}
\multiput(83,0)(20,0){1}{\line(1,0){14}}
\multiput(102.85,-1)(0,2){2}{\line(1,0){14.3}} 
\multiput(59,0)(4,0){6}{\line(1,0){2}} 
\put(110,0){\makebox(0,0){$>$}}
\put(0,-5){\makebox(0,0)[t]{$1$}}
\put(20,-5){\makebox(0,0)[t]{$2$}}
\put(40,-5){\makebox(0,0)[t]{$3$}}
\put(100,-5){\makebox(0,0)[t]{$n\!\! -\!\! 1$}}
\put(120,-5){\makebox(0,0)[t]{$n$}}
\end{picture}
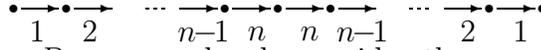
$$
$$
\begin{picture}(180,20)(-10,0)
\multiput( 0,0)(20,0){2}{\circle*{3}}
\multiput(80,0)(20,0){3}{\circle*{3}}
\multiput(180,0)(20,0){2}{\circle*{3}}
\multiput( 3,0)(20,0){2}{\vector(1,0){14}}
\multiput(63,0)(20,0){4}{\vector(1,0){14}}
\multiput(163,0)(20,0){2}{\vector(1,0){14}}
\multiput(50,0)(3,0){3}{\line(1,0){1}} 
\multiput(150,0)(3,0){3}{\line(1,0){1}} 
\put(9,-5){\makebox(0,0)[t]{$1$}}
\put(29,-5){\makebox(0,0)[t]{$2$}}
\put(72,-5){\makebox(0,0)[t]{$n\! \!-\! \!1$}}
\put(92,-7){\makebox(0,0)[t]{$n$}}
\put(112,-7){\makebox(0,0)[t]{$n$}}
\put(132,-5){\makebox(0,0)[t]{$n\! \!-\! \!1$}}
\put(172,-5){\makebox(0,0)[t]{$2$}}
\put(192,-5){\makebox(0,0)[t]{$1$}}
\end{picture}
$$

\eqref{case-Dfin} [$D_n$]
For type $D_n$, we need only consider
 the case that $\aA$ and $\bb = \ee{n-3} \aA$ are \singular~ with 
$\eps(\aA) = \Lambda_{n-3}, \eps(\bb) = \Lambda_{n-2}$.
Let $\cc = \ee{n-2} \bb$. 
By corollaries \ref{cor-parent} and \ref{cor-serre}, 
$\eps(\cc) = m_0 \Lambda_n + m_1 \Lambda_{n-1}$ with
$0 \le m_\ell \le 1$. Let $m = m_0 + m_1$.  
If $m=0$, then $\cc =\v$ so we are done.  
If $m = 1 = m_\ell$, then $\ee{n-\ell} \cc = \v$ and again we are done.
\omitt{
Otherwise, let $\dd = \ee{n} \ee{n-1} \cc = \ee{n-1} \ee{n} \cc.$ 
[fix? ]
If $k \neq n-2$, then $\ee{k} \dd = \ee{n} \ee{n-1} \ee{k} \cc = \zero$.
We can apply Theorem \ref{thm-kash} to see $\ee{n-2}^2 \dd = \zero$,
and so $\dd$ is singular.
}
Otherwise, let $\dd = \ee{n} \ee{n-1} \cc$.
Observe $\epsof{n}(\ee{n-1} \cc) = \epsof{n}(\cc) = 1$
and 
$\epsof{n-1}(\ee{n} \cc) =  1$
 by Lemma \ref{lemma-eps}.
Thus $G(\ee{n} \ee{n-1} \cc) = \eee{n} \eee{n-1} G(\cc) =
\eee{n-1} \eee{n} G(\cc) = G(\ee{n-1} \ee{n} \cc)$, so that
 $ \dd = \ee{n-1} \ee{n} \cc$  as well.
If $k \neq n-2$, then
$\epsof{k}(\dd) = \epsof{k}(\cc) = 0$. 
We can apply Theorem \ref{thm-kash} to see $\ee{n-2}^2 \dd = \zero$,
and so $\dd$ is singular.

 Standard arguments show either $\dd = \v$ or
$\ee{n-2} \dd$ is also singular.  And then this reduces to case
\eqref{case-A}.

\omitt{ types
For type $B_{n}^{(1)}$, we run through the same argument with subscripts
$i$ replaced by $n-i$, and then it reduces to case \eqref{case-Bfin}.
Let $B$ be a highest weight crystal of type $D_n, n \ge 4$, 
$B_{n}^{(1)}, n \ge 3$,
or of  ??? type $A_{2n-1}^{(2) }, n \ge 3$.
}
We again show  the Dynkin diagram and the perfect crystal with
$0$-arrows removed. 
$$ D_n:
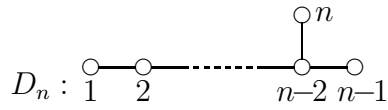
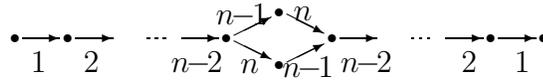
\begin{picture}(106,40)(-5,-10)
\multiput( 0,0)(20,0){2}{\circle{6}}
\multiput(80,0)(20,0){2}{\circle{6}}
\put(80,20){\circle{6}}
\multiput( 3,0)(20,0){2}{\line(1,0){14}}
\multiput(63,0)(20,0){2}{\line(1,0){14}}
\multiput(39,0)(4,0){6}{\line(1,0){2}}
\put(80,3){\line(0,1){14}}
\put(0,-5){\makebox(0,0)[t]{$1$}}
\put(20,-5){\makebox(0,0)[t]{$2$}}
\put(80,-5){\makebox(0,0)[t]{$n\!\! - \!\! 2$}}
\put(103,-5){\makebox(0,0)[t]{$n\!\! -\!\! 1$}}
\put(85,20){\makebox(0,0)[l]{$n$}}
\end{picture}
$$
%
$$
\begin{picture}(215,40)(-20,-15)
\multiput( 0,0)(20,0){2}{\circle*{3}}
\multiput(80,0)(40,0){2}{\circle*{3}}
\multiput(100,-10)(0,20){2}{\circle*{3}}
\multiput(180,0)(20,0){2}{\circle*{3}}
\multiput( 3,0)(20,0){2}{\vector(1,0){14}}
\multiput(63,0)(60,0){2}{\vector(1,0){14}}
\multiput(163,0)(20,0){2}{\vector(1,0){14}}
\multiput(50,0)(3,0){3}{\line(1,0){1}} 
\multiput(150,0)(3,0){3}{\line(1,0){1}} 
\put(9,-5){\makebox(0,0)[t]{$1$}}
\put(29,-5){\makebox(0,0)[t]{$2$}}
\put(69,-5){\makebox(0,0)[t]{$n\! \!-\! \!2$}}
\put(133,-5){\makebox(0,0)[t]{$n\! \!-\! \!2$}}
\put(83,1){\vector(2,1){14}}
\put(83,-1){\vector(2,-1){14}}
\put(103,9){\vector(2,-1){14}}
\put(103,-9){\vector(2,1){14}}
\put(172,-5){\makebox(0,0)[t]{$2$}}
\put(192,-5){\makebox(0,0)[t]{$1$}}
\put(86,13){\makebox(0,0)[t]{$n\! \!-\! \!1$}}
\put(89,-8){\makebox(0,0)[t]{$n$}}
\put(111,-7){\makebox(0,0)[t]{$n\! \!-\! \!1$}}
\put(109,12){\makebox(0,0)[t]{$n$}}
\end{picture}
$$

\omitt{
\eqref{case-D2}
Everything is covered in the cases \ref{case-D} and \ref{case-C} above.
\omitt{ types
$D_{n}^{(1)}, n \ge 4$ or  $A_{2n-1}^{(2)}, n \ge 3$.
}
We list the Dynkin diagrams below
$$ D_n^{(1)}:
\begin{picture}(106,40)(-5,-5)
\multiput( 0,0)(20,0){2}{\circle{6}}
\multiput(80,0)(20,0){2}{\circle{6}}
\multiput(20,20)(60,0){2}{\circle{6}}
\multiput( 3,0)(20,0){2}{\line(1,0){14}}
\multiput(63,0)(20,0){2}{\line(1,0){14}}
\multiput(39,0)(4,0){6}{\line(1,0){2}}
\multiput(20,3)(60,0){2}{\line(0,1){14}}
\put(0,-5){\makebox(0,0)[t]{$1$}}
\put(20,-5){\makebox(0,0)[t]{$2$}}
\put(80,-5){\makebox(0,0)[t]{$n\!\! - \!\! 2$}}
\put(103,-5){\makebox(0,0)[t]{$n\!\! -\!\! 1$}}
\put(25,20){\makebox(0,0)[l]{$0$}}
\put(85,20){\makebox(0,0)[l]{$n$}}
\end{picture}
\qquad
A_{2n-1}^{(2)}:
\begin{picture}(126,40)(-5,-5)
\multiput( 0,0)(20,0){3}{\circle{6}}
\multiput(100,0)(20,0){2}{\circle{6}}
\put(20,20){\circle{6}}
\multiput( 3,0)(20,0){3}{\line(1,0){14}}
\multiput(83,0)(20,0){1}{\line(1,0){14}}
\put(20,3){\line(0,1){14}}
\multiput(102.85,-1)(0,2){2}{\line(1,0){14.3}} 
\multiput(59,0)(4,0){6}{\line(1,0){2}} 
\put(110,0){\makebox(0,0){$<$}}
\put(0,-5){\makebox(0,0)[t]{$1$}}
\put(20,-5){\makebox(0,0)[t]{$2$}}
\put(40,-5){\makebox(0,0)[t]{$3$}}
\put(100,-5){\makebox(0,0)[t]{$n\!\! -\!\! 1$}}
\put(120,-7){\makebox(0,0)[t]{$n$}}
\put(25,20){\makebox(0,0)[l]{$0$}}
\end{picture} \; .
$$
}
\eqref{case-therest.unique}
[ $C_n^{(1)},
A_{2n}^{(2)},
A_{2n}^{(2) \dagger},
D_{n+1}^{(2)}$]
The local nature of \singular ity allows us to apply the results
from cases \eqref{case-A}, \eqref{case-C}, \eqref{case-Bfin}
to these types (sometimes reindexing $i$ for $n-i$ when encountering
$\Lambda_0$). 

We list the perfect crystals, and for completeness, the possible
 \walks~ from  $\v$ to $\aA$.

$C_n^{(1)}, (n \ge 2)$
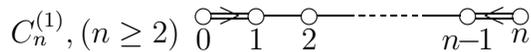
\begin{picture}(126,25)(-5,-10)
\multiput( 0,0)(20,0){3}{\circle{6}}
\multiput(100,0)(20,0){2}{\circle{6}}
\multiput(23,0)(20,0){2}{\line(1,0){14}}
\put(83,0){\line(1,0){14}}
\multiput( 2.85,-1)(0,2){2}{\line(1,0){14.3}} 
\multiput(102.85,-1)(0,2){2}{\line(1,0){14.3}} 
\multiput(59,0)(4,0){6}{\line(1,0){2}} 
\put(10,0){\makebox(0,0){$>$}} \put(110,0){\makebox(0,0){$<$}}
\put(0,-5){\makebox(0,0)[t]{$0$}}
\put(20,-5){\makebox(0,0)[t]{$1$}}
\put(40,-5){\makebox(0,0)[t]{$2$}}
\put(100,-5){\makebox(0,0)[t]{$n\!\! -\!\! 1$}}
\put(120,-5){\makebox(0,0)[t]{$n$}}
\end{picture}

\begin{picture}(200,52)(-10,-15)
\multiput( 0,0)(20,0){2}{\circle*{3}}
\multiput(80,0)(20,0){2}{\circle*{3}}
\multiput(160,0)(20,0){2}{\circle*{3}}
\multiput( 3,0)(20,0){2}{\vector(1,0){14}}
\multiput(63,0)(20,0){3}{\vector(1,0){14}}
\multiput(143,0)(20,0){2}{\vector(1,0){14}}
\multiput(50,0)(3,0){3}{\line(1,0){1}} 
\multiput(130,0)(3,0){3}{\line(1,0){1}} 
\put(91,4){\oval(178,25)[t]}
\put(99,16){\vector(-1,0){3}}
\put(9,-5){\makebox(0,0)[t]{$1$}}
\put(29,-5){\makebox(0,0)[t]{$2$}}
\put(72,-5){\makebox(0,0)[t]{$n\! \!-\! \!1$}}
\put(92,-7){\makebox(0,0)[t]{$n$}}
\put(112,-5){\makebox(0,0)[t]{$n\! \!-\! \!1$}}
\put(152,-5){\makebox(0,0)[t]{$2$}}
\put(172,-5){\makebox(0,0)[t]{$1$}}
\put(105,25){\makebox(0,0)[t]{$0$}}
\end{picture}

\noindent
$\aA = \ffi \ff{i + 1}  \cdots \ff{n -1} \ff{n} \ff{n-1}
 \cdots \ff{1} \ff{0} \ff{1} \cdots \ff{n} \cdots \ff{0} \cdots  \v$
\newline \noindent
$\aA = \ffi \ff{i - 1}  \cdots \ff{1} \ff{0} \ff{1} \cdots \ff{n} \cdots
\ff{0} \cdots  \ff{n}  \cdots \v$

\vspace{5ex}
%
%
$A_{2n}^{(2)}, (n \ge 2)$
\begin{picture}(126,25)(-5,-15)
\multiput( 0,0)(20,0){3}{\circle{6}}
\multiput(100,0)(20,0){2}{\circle{6}}
\multiput(23,0)(20,0){2}{\line(1,0){14}}
\put(83,0){\line(1,0){14}}
\multiput( 2.85,-1)(0,2){2}{\line(1,0){14.3}} 
\multiput(102.85,-1)(0,2){2}{\line(1,0){14.3}} 
\multiput(59,0)(4,0){6}{\line(1,0){2}} 
\put(10,0){\makebox(0,0){$<$}} \put(110,0){\makebox(0,0){$<$}}
\put(0,-5){\makebox(0,0)[t]{$0$}}
\put(20,-5){\makebox(0,0)[t]{$1$}}
\put(40,-5){\makebox(0,0)[t]{$2$}}
\put(100,-5){\makebox(0,0)[t]{$n\!\! -\!\! 1$}}
\put(120,-5){\makebox(0,0)[t]{$n$}}
\end{picture}

\begin{picture}(200,55)(-10,-20)
\multiput( 0,0)(20,0){2}{\circle*{3}}
\multiput(80,0)(20,0){2}{\circle*{3}}
\multiput(160,0)(20,0){2}{\circle*{3}}
\multiput( 3,0)(20,0){2}{\vector(1,0){14}}
\multiput(63,0)(20,0){3}{\vector(1,0){14}}
\multiput(143,0)(20,0){2}{\vector(1,0){14}}
\multiput(50,0)(3,0){3}{\line(1,0){1}} 
\multiput(130,0)(3,0){3}{\line(1,0){1}} 
\put(90,23){\circle*{3}}
\put(2,2){\vector(4,1){84}}
\put(93,24){\vector(4,-1){84}}
\put(9,-5){\makebox(0,0)[t]{$1$}}
\put(29,-5){\makebox(0,0)[t]{$2$}}
\put(72,-5){\makebox(0,0)[t]{$n\! \!-\! \!1$}}
\put(92,-7){\makebox(0,0)[t]{$n$}}
\put(112,-5){\makebox(0,0)[t]{$n\! \!-\! \!1$}}
\put(152,-5){\makebox(0,0)[t]{$2$}}
\put(172,-5){\makebox(0,0)[t]{$1$}}
\put(135,23){\makebox(0,0)[t]{$0$}}
\put(45,23){\makebox(0,0)[t]{$0$}}
\end{picture}

\noindent
$\aA = \ffi \ff{i + 1}  \cdots \ff{n-2}  \ff{n-1} \ff{ n} \ff{n-1}  \ff{n-2}
\cdots  \ff{2}
\ff{1} \ff{0} \ff{0} \ff{1} \ff{2}  \cdots
\ff{n-1} \ff{ n} \ff{n-1} \cdots \ff{1} \ff{0}$   $\ff{0} \ff{1}
\cdots \ff{n-1} \ff{ n} \ff{n-1} \cdots \v$
\newline \noindent
$\aA = \ffi \ff{i - 1}  \cdots
\cdots  \ff{2} \ff{1} \ff{0} \ff{0} \ff{1} \ff{2}  \cdots
\ff{n-1} \ff{ n} \ff{n-1} \cdots \ff{1} \ff{0} \ff{0} \ff{1}
\cdots \v$
\vspace{5ex}

$A_{2n}^{(2)\dagger}, (n \ge 2)$
\begin{picture}(126,25)(-5,-15)
\multiput( 0,0)(20,0){3}{\circle{6}}
\multiput(100,0)(20,0){2}{\circle{6}}
\multiput(23,0)(20,0){2}{\line(1,0){14}}
\put(83,0){\line(1,0){14}}
\multiput( 2.85,-1)(0,2){2}{\line(1,0){14.3}} 
\multiput(102.85,-1)(0,2){2}{\line(1,0){14.3}} 
\multiput(59,0)(4,0){6}{\line(1,0){2}} 
\put(10,0){\makebox(0,0){$>$}} \put(110,0){\makebox(0,0){$>$}}
\put(0,-5){\makebox(0,0)[t]{$0$}}
\put(20,-5){\makebox(0,0)[t]{$1$}}
\put(40,-5){\makebox(0,0)[t]{$2$}}
\put(100,-5){\makebox(0,0)[t]{$n\!\! -\!\! 1$}}
\put(120,-5){\makebox(0,0)[t]{$n$}}
\put(0,13){\makebox(0,0)[t]{$$}}
\put(20,13){\makebox(0,0)[t]{$$}}
\put(40,13){\makebox(0,0)[t]{$$}}
\put(100,13){\makebox(0,0)[t]{$$}}
\put(120,13){\makebox(0,0)[t]{$$}}
\put(120,13){\makebox(0,0)[t]{$$}}
\end{picture}

%
\begin{picture}(200,45)(-10,-15)
\multiput( 0,0)(20,0){2}{\circle*{3}}
\multiput(80,0)(20,0){3}{\circle*{3}}
\multiput(180,0)(20,0){2}{\circle*{3}}
\multiput( 3,0)(20,0){2}{\vector(1,0){14}}
\multiput(63,0)(20,0){4}{\vector(1,0){14}}
\multiput(163,0)(20,0){2}{\vector(1,0){14}}
\multiput(50,0)(3,0){3}{\line(1,0){1}} 
\multiput(150,0)(3,0){3}{\line(1,0){1}} 
\put(101,4){\oval(198,25)[t]}
\put(99,16){\vector(-1,0){3}}
\put(9,-5){\makebox(0,0)[t]{$1$}}
\put(29,-5){\makebox(0,0)[t]{$2$}}
\put(72,-5){\makebox(0,0)[t]{$n\! \!-\! \!1$}}
\put(92,-7){\makebox(0,0)[t]{$n$}}
\put(112,-7){\makebox(0,0)[t]{$n$}}
\put(132,-5){\makebox(0,0)[t]{$n\! \!-\! \!1$}}
\put(172,-5){\makebox(0,0)[t]{$2$}}
\put(192,-5){\makebox(0,0)[t]{$1$}}
\put(105,25){\makebox(0,0)[t]{$0$}}
\end{picture}

$\aA = \ffi \ff{i + 1}  \cdots \ff{n-2} 
\ff{n-1} \ff{ n}  \ff{n} \ff{n-1} 
\ff{n-2} \cdots  \ff{2} \ff{1}  \ff{0} \ff{1} \ff{2}  \cdots
\ff{n-1} \ff{ n}  \ff{n} \ff{n-1} 
\cdots$  $ \ff{1}  \ff{0} \ff{1} \cdots \v$
\newline \noindent
$\aA = \ffi \ff{i - 1}  \cdots
 \ff{2} \ff{1}  \ff{0} \ff{1} \ff{2}  \cdots \ff{n-1} \ff{ n}  \ff{n} \ff{n-1} 
\cdots \ff{1}  \ff{0} \ff{1} \cdots 
\ff{n-1} \ff{ n}  \ff{n} \ff{n-1} \cdots \v.$

\vspace{5ex}
$D_{n+1}^{(2)}, (n \ge 2)$
\begin{picture}(126,20)(-5,-5)
\multiput( 0,0)(20,0){3}{\circle{6}}
\multiput(100,0)(20,0){2}{\circle{6}}
\multiput(23,0)(20,0){2}{\line(1,0){14}}
\put(83,0){\line(1,0){14}}
\multiput( 2.85,-1)(0,2){2}{\line(1,0){14.3}} 
\multiput(102.85,-1)(0,2){2}{\line(1,0){14.3}} 
\multiput(59,0)(4,0){6}{\line(1,0){2}} 
\put(10,0){\makebox(0,0){$<$}} \put(110,0){\makebox(0,0){$>$}}
\put(0,-5){\makebox(0,0)[t]{$0$}}
\put(20,-5){\makebox(0,0)[t]{$1$}}
\put(40,-5){\makebox(0,0)[t]{$2$}}
\put(100,-5){\makebox(0,0)[t]{$n\!\! -\!\! 1$}}
\put(120,-5){\makebox(0,0)[t]{$n$}}
\end{picture}

\begin{picture}(190,52)(-10,-15)
\multiput( 0,0)(20,0){2}{\circle*{3}}
\multiput(80,0)(20,0){3}{\circle*{3}}
\multiput(180,0)(20,0){2}{\circle*{3}}
\multiput( 3,0)(20,0){2}{\vector(1,0){14}}
\multiput(63,0)(20,0){4}{\vector(1,0){14}}
\multiput(163,0)(20,0){2}{\vector(1,0){14}}
\multiput(50,0)(3,0){3}{\line(1,0){1}} 
\multiput(150,0)(3,0){3}{\line(1,0){1}} 
\put(9,-5){\makebox(0,0)[t]{$1$}}
\put(29,-5){\makebox(0,0)[t]{$2$}}
\put(72,-5){\makebox(0,0)[t]{$n\! \!-\! \!1$}}
\put(92,-7){\makebox(0,0)[t]{$n$}}
\put(112,-7){\makebox(0,0)[t]{$n$}}
\put(132,-5){\makebox(0,0)[t]{$n\! \!-\! \!1$}}
\put(172,-5){\makebox(0,0)[t]{$2$}}
\put(192,-5){\makebox(0,0)[t]{$1$}}
\put(100,26){\circle*{3}}
\put(2,2){\vector(4,1){93}}
\put(103,26){\vector(4,-1){93}}
\put(150,23){\makebox(0,0)[t]{$0$}}
\put(45,23){\makebox(0,0)[t]{$0$}}
\end{picture}

\noindent
$\aA = \ffi \ff{i + 1}  \cdots \ff{n-1} \ff{n} \ff{n} \ff{n-1} 
\cdots \ff{2} \ff{1}  \ff{0} \ff{0} \ff{1} \ff{2}  \cdots
\ff{ n}  \ff{n}  \cdots \v$
\newline \noindent
$\aA = \ffi \ff{i - 1}  \cdots
 \ff{2} \ff{1}  \ff{0} \ff{0} \ff{1} \ff{2}  \cdots
\ff{n-1} \ff{n} \ff{n} \ff{n-1} \cdots \ff{0} \ff{0} \cdots \v$

\vspace{5ex}
\eqref{case-therest.infinite}
[$D_{n}^{(1)}, 
A_{2n-1}^{(2)},
B_{n}^{(1)}$]
As above, we may apply the results
from cases \eqref{case-A}, \eqref{case-C}, \eqref{case-Bfin}, \eqref{case-Dfin}
to these types (with appropriate reindexing).

We  list the perfect crystals, and for completeness, the possible
 \walks~ from $\v$ to $\aA$.

Below we again  use the notation
$\ff{n-1, n}$ to stand for  either
$ \ff{n-1} \ff{n} $ or $ \ff{n} \ff{n-1} $  in types
$D_{n},  D_{n}^{(1)}$ and
$\ff{0, 1}$ to stand for  either
$ \ff{0} \ff{1} $ or $ \ff{1} \ff{0}$ in types  $ D_{n}^{(1)}, B_{n}^{(1)}$.

$D_n^{(1)}, (n \ge 4)$  
\begin{picture}(106,40)(-5,-5)
\multiput( 0,0)(20,0){2}{\circle{6}}
\multiput(80,0)(20,0){2}{\circle{6}}
\multiput(20,20)(60,0){2}{\circle{6}} \multiput(
3,0)(20,0){2}{\line(1,0){14}}
\multiput(63,0)(20,0){2}{\line(1,0){14}}
\multiput(39,0)(4,0){6}{\line(1,0){2}}
\multiput(20,3)(60,0){2}{\line(0,1){14}}
\put(0,-5){\makebox(0,0)[t]{$1$}}
\put(20,-5){\makebox(0,0)[t]{$2$}}
\put(80,-5){\makebox(0,0)[t]{$n\!\! - \!\! 2$}}
\put(103,-5){\makebox(0,0)[t]{$n\!\! -\!\! 1$}}
\put(25,20){\makebox(0,0)[l]{$0$}}
\put(85,20){\makebox(0,0)[l]{$n$}}
\end{picture}

\begin{picture}(215,65)(-10,-30)
\multiput( 0,0)(20,0){2}{\circle*{3}}
\multiput(80,0)(40,0){2}{\circle*{3}}
\multiput(100,-10)(0,20){2}{\circle*{3}}
\multiput(180,0)(20,0){2}{\circle*{3}}
\multiput( 3,0)(20,0){2}{\vector(1,0){14}}
\multiput(63,0)(60,0){2}{\vector(1,0){14}}
\multiput(163,0)(20,0){2}{\vector(1,0){14}}
\multiput(50,0)(3,0){3}{\line(1,0){1}} 
\multiput(150,0)(3,0){3}{\line(1,0){1}} 
\put(91,-4){\oval(178,25)[b]}
\put(111,4){\oval(178,25)[t]}
\put(149,16){\vector(-1,0){3}}
\put(49,-17){\vector(-1,0){3}}
\put(9,-5){\makebox(0,0)[t]{$1$}}
\put(29,-5){\makebox(0,0)[t]{$2$}}
\put(69,-5){\makebox(0,0)[t]{$n\! \!-\! \!2$}}
\put(133,-5){\makebox(0,0)[t]{$n\! \!-\! \!2$}}
\put(83,1){\vector(2,1){14}}
\put(83,-1){\vector(2,-1){14}}
\put(103,9){\vector(2,-1){14}}
\put(103,-9){\vector(2,1){14}}
\put(172,-5){\makebox(0,0)[t]{$2$}}
\put(192,-5){\makebox(0,0)[t]{$1$}}
\put(155,25){\makebox(0,0)[t]{$0$}}
\put(55,-22){\makebox(0,0)[t]{$0$}}
\put(86,13){\makebox(0,0)[t]{$n\! \!-\! \!1$}}
\put(89,-8){\makebox(0,0)[t]{$n$}}
\put(111,-7){\makebox(0,0)[t]{$n\! \!-\! \!1$}}
\put(109,12){\makebox(0,0)[t]{$n$}}
\end{picture}

\noindent
$\aA = \ffi \ff{i + 1}  \cdots \ff{n-2} \ff{n,n-1} \ff{n-2}
\cdots \ff{2} \ff{0,1} \ff{2}  \cdots  \ff{n,n-1} \cdots \v$
\newline \noindent
$\aA = \ffi \ff{i - 1}  \cdots \ff{2}
\ff{1}  \ff{0} \ff{0} \ff{1}
 \ff{2}  \cdots \ff{n-1} \ff{n} \ff{n} \ff{n-1} 
\cdots \ff{0} \ff{0} \cdots \v$

\vspace{5ex}

$A_{2n}^{(2)}, (n \ge 2)$  
\begin{picture}(126,20)(-5,-5)
\multiput( 0,0)(20,0){3}{\circle{6}}
\multiput(100,0)(20,0){2}{\circle{6}}
\multiput(23,0)(20,0){2}{\line(1,0){14}}
\put(83,0){\line(1,0){14}}
\multiput( 2.85,-1)(0,2){2}{\line(1,0){14.3}} 
\multiput(102.85,-1)(0,2){2}{\line(1,0){14.3}} 
\multiput(59,0)(4,0){6}{\line(1,0){2}} 
\put(10,0){\makebox(0,0){$<$}} \put(110,0){\makebox(0,0){$<$}}
\put(0,-5){\makebox(0,0)[t]{$0$}}
\put(20,-5){\makebox(0,0)[t]{$1$}}
\put(40,-5){\makebox(0,0)[t]{$2$}}
\put(100,-5){\makebox(0,0)[t]{$n\!\! -\!\! 1$}}
\put(120,-5){\makebox(0,0)[t]{$n$}}
\put(0,13){\makebox(0,0)[t]{$2$}}
\put(20,13){\makebox(0,0)[t]{$2$}}
\put(40,13){\makebox(0,0)[t]{$2$}}
\put(100,13){\makebox(0,0)[t]{$2$}}
\put(120,13){\makebox(0,0)[t]{$2$}}
\put(120,13){\makebox(0,0)[t]{$2$}}
\end{picture}

\begin{picture}(210,65)(-7,-30)
\multiput( 0,0)(20,0){2}{\circle*{3}}
\multiput(90,0)(20,0){2}{\circle*{3}}
\multiput(180,0)(20,0){2}{\circle*{3}}
\multiput( 3,0)(20,0){2}{\vector(1,0){14}}
\multiput(73,0)(20,0){4}{\vector(1,0){14}}
\multiput(163,0)(20,0){2}{\vector(1,0){14}}
\multiput(50,0)(3,0){3}{\line(1,0){1}} 
\multiput(150,0)(3,0){3}{\line(1,0){1}} 
\put(91,-4){\oval(178,25)[b]}
\put(111,4){\oval(178,25)[t]}
\put(149,16){\vector(-1,0){3}}
\put(49,-17){\vector(-1,0){3}}
\put(9,-5){\makebox(0,0)[t]{$1$}}
\put(29,-5){\makebox(0,0)[t]{$2$}}
\put(80,-5){\makebox(0,0)[t]{$n\! \!-\! \!1$}}
\put(100,-7){\makebox(0,0)[t]{$n$}}
\put(120,-5){\makebox(0,0)[t]{$n\! \!-\! \!1$}}
\put(172,-5){\makebox(0,0)[t]{$2$}}
\put(192,-5){\makebox(0,0)[t]{$1$}}
\put(155,25){\makebox(0,0)[t]{$0$}}
\put(55,-22){\makebox(0,0)[t]{$0$}}
\end{picture}

\noindent
$\aA = \ffi \ff{i + 1}  \cdots \ff{n-2}  
\ff{n-1} \ff{ n} \ff{n-1} \ff{n-2}
\cdots  \ff{2} \ff{0, 1}  \ff{2}  \cdots
\ff{n-1} \ff{ n} \ff{n-1} \cdots \ff{0, 1} \cdots \v$
\newline \noindent
$\aA = \ffi \ff{i - 1}  \cdots
 \ff{2} \ff{0, 1}  \ff{2}  \cdots
\ff{n-1} \ff{ n} \ff{n-1} 
\cdots \ff{0, 1} \cdots \v$
%

$B_n^{(1)}, (n \ge 3)$
\begin{picture}(126,40)(-5,-5)
\multiput( 0,0)(20,0){3}{\circle{6}}
\multiput(100,0)(20,0){2}{\circle{6}} \put(20,20){\circle{6}}
\multiput( 3,0)(20,0){3}{\line(1,0){14}}
\multiput(83,0)(20,0){1}{\line(1,0){14}}
\put(20,3){\line(0,1){14}}
\multiput(102.85,-1)(0,2){2}{\line(1,0){14.3}} 
\multiput(59,0)(4,0){6}{\line(1,0){2}} 
\put(110,0){\makebox(0,0){$>$}} \put(0,-5){\makebox(0,0)[t]{$1$}}
\put(20,-5){\makebox(0,0)[t]{$2$}}
\put(40,-5){\makebox(0,0)[t]{$3$}}
\put(100,-5){\makebox(0,0)[t]{$n\!\! -\!\! 1$}}
\put(120,-5){\makebox(0,0)[t]{$n$}}
\put(25,20){\makebox(0,0)[l]{$0$}}
\put(120,13){\makebox(0,0)[t]{$2$}}
\end{picture}

\begin{picture}(180,60)(-10,-30)
\multiput( 0,0)(20,0){2}{\circle*{3}}
\multiput(80,0)(20,0){3}{\circle*{3}}
\multiput(180,0)(20,0){2}{\circle*{3}}
\multiput( 3,0)(20,0){2}{\vector(1,0){14}}
\multiput(63,0)(20,0){4}{\vector(1,0){14}}
\multiput(163,0)(20,0){2}{\vector(1,0){14}}
\multiput(50,0)(3,0){3}{\line(1,0){1}} 
\multiput(150,0)(3,0){3}{\line(1,0){1}} 
\put(91,-4){\oval(178,25)[b]}
\put(111,4){\oval(178,25)[t]}
\put(149,16){\vector(-1,0){3}}
\put(49,-17){\vector(-1,0){3}}
\put(9,-5){\makebox(0,0)[t]{$1$}}
\put(29,-5){\makebox(0,0)[t]{$2$}}
\put(72,-5){\makebox(0,0)[t]{$n\! \!-\! \!1$}}
\put(92,-7){\makebox(0,0)[t]{$n$}}
\put(112,-7){\makebox(0,0)[t]{$n$}}
\put(132,-5){\makebox(0,0)[t]{$n\! \!-\! \!1$}}
\put(172,-5){\makebox(0,0)[t]{$2$}}
\put(192,-5){\makebox(0,0)[t]{$1$}}
\put(155,25){\makebox(0,0)[t]{$0$}}
\put(55,-22){\makebox(0,0)[t]{$0$}}
\end{picture}
%

\noindent
$\aA = \ffi \ff{i + 1}  \cdots \ff{n-2}  \ff{n-1} \ff{ n} \ff{ n} \ff{n-1} 
 \ff{n-2} \cdots  \ff{2} \ff{0, 1}  \ff{2}  \cdots
\ff{n-1} \ff{ n} \ff{ n} \ff{n-1} \cdots \ff{0, 1} \cdots \v$
\newline \noindent
$\aA = \ffi \ff{i - 1}  \cdots
 \ff{2} \ff{0, 1}  \ff{2}  \cdots \ff{n-1} \ff{ n} \ff{n-1}  \cdots \v$

\end{proof}
\section{Existence}
\label{sec-existence}
The above theorems consisted of several ``uniqueness" statements.
The corresponding existence statements also hold. 

In Theorem \ref{thm-global}, we described all
sequences $i_1, i_2, \cdots, i_k$
such that 
 $ \aA =  \ff{i_1} \ff{i_2} \cdots \ff{i_k} \v,$
where $\v$ is the highest weight node and we required $\aA$ to
be \singular~ with \singular~ parent. 
These possible sequences corresponded to
\walks~ $\xrightarrow{i_1} \xrightarrow{i_2}\cdots \xrightarrow{i_k}$
on a perfect crystal.
Below we will exhibit a highest weight crystal (one of level $1$
or level $2$ suffices) and such a node $\aA$ for every such
\walk~ (excluding of course \walks~ where $a_{i_1, i_2} \ge 0$, as
in that case $\aA$ would not be \singular).

We recall that the tensor product of crystals 
$B_2 \otimes B_1$ is defined by the nodes being the Cartesian
product of the nodes of $B_2$ and $B_1$,  $\wt(\bb_2 \otimes \bb_1)
= \wt(\bb_2) + \wt(\bb_1)$, and arrows are described by the
following rule 
$$\ei(b_2 \otimes b_1) = \begin{cases}
	\ei b_2 \otimes b_1 & \text{ if $\varphi_i(b_2) \ge \epsi(b_1)$} \\
	b_2 \otimes \ei b_1 & \text{ if $\varphi_i(b_2) < \epsi(b_1)$.}
 \end{cases}$$
Consequently
\begin{gather}
\label{eqn-eps}
\epsi(b_2 \otimes b_1) = \epsi(b_2) + \max\{0, \epsi(b_1) - \varphi_i(b_2) \} \\
\varphi_i(b_2 \otimes b_1)
	= \varphi_i(b_1) + \max\{0, \varphi_i(b_2) - \epsi(b_1) \}.
\end{gather}
We recall the following theorem.
\begin{theorem}[\cite{KK.core}, \cite{KKMMNN.perfect} ]
\label{thm-KMN}
Let $\lambda$ be a dominant integral weight of level $k$, and $B$
be a perfect crystal of level $\ell$, and suppose $k \ge \ell$.
Then
$$B(\lambda) \otimes B \simeq
	\bigoplus_{\bb \in B^{\le \lambda}}
B(\lambda + \wt(\bb) )$$
where $B^{\le \lambda} = \{ \bb \in B \mid \epsi(\bb) \le \langle h_i,
\lambda \rangle \, \forall i \}$.
\end{theorem}
We set
$$\psi_k^{\lambda, \mu} : B(\mu) \to B(\lambda) \otimes (B^{1,1})^{\otimes k}
$$
to be the embedding dictated by the above theorem, when it is defined.
Observe
$(\psi_{k'}^{\lambda, \nu} \otimes {id}^{\otimes k}) \circ
\psi_k^{\nu, \mu} = \psi_{k' + k}^{\lambda + \mu}$.

We refer the reader to the appendix for a list of the level $1$
perfect crystals $B^{1,1}$ (including $B^{n,1}$ in type $A$).
There is a standard way of labelling the nodes, but it will be
convenient here to ignore that convention, so we have ommitted
that labelling in the appendix.

Let 
$\xrightarrow{i_1} \xrightarrow{i_2}\cdots \xrightarrow{i_k}$
be a \walk~ on $B^{1,1}$ (or $B^{n,1}$ in type $A$).
Let
\begin{gather}
m = | \{ r \mid a_{i_r, i_{r+1}} \ge 0, 1 \le r < k \} |
\end{gather}
and let
\begin{gather*}
\jboxi{1} \* \jboxi{2} \* \cdots \* \jboxi{k} 
\in (B^{1,1})^{\* (k-m) }
\end{gather*}
be such that $\jboxi{r}$ is the node
$$ \xrightarrow{i_r} \jboxi{r} \xrightarrow{i_{r+1}}$$
 with an $i_r$-colored
arrow going in and $i_{r+1}$ going out {\it if\/} $a_{i_r, i_{r+1}} < 0$.
So long as $k > 1$, these nodes are well-defined and this also
determines $\jboxi{k}$.
Observe that in the case $a_{i_1, i_2} \ge 0$, the node we
describe is actually then 
$\jboxi{2} \*  \cdots \* \jboxi{k}$.
Also note the labelling very much depends on the \walk, and that
one node can  receive many different labels.


\begin{lemma}
\label{lemma-eps-tensor}
Let 
$\jboxi{1} \* \jboxi{2} \* \cdots \* \jboxi{k} 
\in (B^{1,1})^{\* (k-m) }$ be as above. For all $i \in I$,
\begin{enumerate}
\item
$\epsi(
\jboxi{1} \* \jboxi{2} \* \cdots \* \jboxi{k} )
=  \epsi( \jboxi{1})$ 
\item
$\ei(\jboxi{1} \* \jboxi{2} \* \cdots \* \jboxi{k} )
= \ei (\jboxi{1}) \* \jboxi{2} \* \cdots \* \jboxi{k}$ 
\item
$\varphi_i(\jboxi{1} \* \jboxi{2} \* \cdots \* \jboxi{k} )
=  \varphi_i( \jboxi{k})$ 
\end{enumerate}
\end{lemma}
\begin{proof}
When $k=1$ this is immediate.
Recall from \eqref{eqn-eps}, 
$\epsi( \jboxi{1} \* \jboxi{2} \* \cdots \* \jboxi{k} )
 = \epsi(\jboxi{1}) + \max\{0, \epsi(\jboxi{2} \* \cdots \* \jboxi{k} )
- \varphi_i(\jboxi{1}) \}.$
By the inductive hypothesis, 
$\epsi(  \jboxi{2} \* \cdots \* \jboxi{k} )
=  \epsi( \jboxi{2})$ (by which we mean the leftmost node in case $a_{i_2, i_3}
\ge 0$).
If $ \epsi( \jboxi{2}) = 0$, we are done.
If $ \epsi( \jboxi{2}) \neq 0$, we will show $ \epsi( \jboxi{2})
-\varphi_i(\jboxi{1}) \le 0$.

Consider the following possibilities.
First, $i = i_2$ and 
$\epsof{i_2}(\jboxi{2}) = 1$.
As $\varphi_{i_2}(\jboxi{1}) \ge 1$, we are done.
Second, suppose
 $i = i_2$ and 
 $\epsof{i_2}(\jboxi{2}) > 1.$  
In fact, because we assume $\jboxi{1}$ contributes to the tensor and 
$\xrightarrow{i_2}$ joins $\jboxi{1}$ to $\jboxi{2}$, this cannot happen.
It would mean $\jboxi{2}$ does not contribute,
and the ``leftmost" node we refer to above is in fact $\jboxi{3}$.
We have 
$$ \xrightarrow{i_1} \jboxi{1}  \xrightarrow{i_2} \cdot
	\xrightarrow{i_3 = i_2} \jboxi{3},$$
as $a_{i_2, i_3} = 2 \ge 0$,
so that  
$\epsof{i_2}(\widehat{ \jboxi{2} } \* \jboxi{3} \* \cdots \jboxi{k}) =
\epsof{i_2}(\jboxi{3}) = \varphi_{i_2}(\jboxi{1}) = 2$.

Third, suppose $i \neq i_2$.  Then we must have 
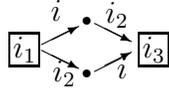
\begin{center}
\begin{picture}(215,35)(-10,-14)
\multiput(100,-10)(0,20){2}{\circle*{3}}
\put(83,1){\vector(2,1){14}}
\put(83,-1){\vector(2,-1){14}}
\put(103,9){\vector(2,-1){14}}
\put(103,-9){\vector(2,1){14}}
\put(86,13){\makebox(5,5)[t]{$i$}}
\put(89,-10){\makebox(5,5)[t]{$i_2$}}
\put(111,-9){\makebox(5,5)[t]{$i$}}
\put(109,12){\makebox(5,5)[t]{$i_2$}}
\put(71,-7){\framebox(11,12)[c]{$i_1$}}
\put(120,-7){\framebox(11,12)[c]{$i_3$}}
\end{picture}
\end{center}
and 
$ \epsof{i}(\jboxi{3}) = \varphi_{i}(\jboxi{1}) = 1$.
Again, $\jboxi{2}$ does not contribute.

Computing $\varphi_i$ is similar.
The above conclusions come from examining all $B^{1,1}$ and from our
definition of the node in $(B^{1,1})^{\otimes (k-m)}$ that our
\walk~ specifies.

The rule for computing $\ei$ of a tensor product gives us the second
statement.
\end{proof}

Write $\v_\lambda \in B(\lambda)$ for the highest weight node.
\begin{proposition}
\label{prop-exist}
Let 
$\xrightarrow{i_1} \xrightarrow{i_2}\cdots \xrightarrow{i_k}$
be a \walk~ on $B^{1,1}$ (or $B^{n,1}$ in type $A$).
Let $\jboxi{0}$ be the node such that $\jboxi{0} \xrightarrow{i_1}
\jboxi{1}$, and let $\lambda = \varepsilon(\jboxi{0})$.
Let $\mu = \varphi(\jboxi{k\! -\! 1})$.
\begin{enumerate}
\item
$\v_\lambda \* \jboxi{0} \* \jboxi{1} \* \cdots \* \jboxi{k\!-\! 1}
= \psi_{k-m}^{\lambda,\mu}( \v_\mu)$
\item
$\v_\lambda \* \jboxi{1} \*  \cdots \* \jboxi{k}
= \psi_{k-m}^{\lambda,\mu}( \aA)$,
where $\aA = 
    \ff{i_1} \ff{i_2} \cdots \ff{i_k} \v_\mu.$
\item
In particular  $\aA \neq \zero$, and if $ \jboxi{1} \*  \cdots \* \jboxi{k}$
is \singular~ (with \singular~ parent) so is $\aA$.
\end{enumerate}
\end{proposition}

\begin{proof}
1. From \eqref{eqn-eps},
$\epsi(\v_\lambda \* \jboxi{0} \* \jboxi{1} \* \cdots \* \jboxi{k\! -\! 1})
=  \epsi(\v_\lambda) + \max\{0,
 \epsi(\jboxi{0} \*  \cdots \* \jboxi{k\! -\! 1})
- \varphi_i(\v_\lambda) \}
=  0 + \max\{0, \epsi(\jboxi{0}) - \varphi_i(\v_\lambda) \}
=0$
for all $i$
by our choice of $\lambda$.
Hence it is a highest weight node.  
Lemma \ref{lemma-eps-tensor} computes its weight is $\mu$,
so it must be the image of $\v_\mu$.
Notice $\mu$ is of level $1$ or $2$.

2. We only need show
$\v_\lambda \* \jboxi{1} \*  \cdots \* \jboxi{k}
=
\ff{i_1} \ff{i_2} \cdots \ff{i_k} 
(\v_\lambda \* \jboxi{0} \* \jboxi{1} \* \cdots \* \jboxi{k\! -\! 1})$.
We will induct on $k$.

In the case $k=1$, $\ei (\v_\lambda \* \jboxi{1} ) = \zero$ if $i \neq i_1$
and $\ee{i_1}(\v_\lambda \* \jboxi{1} )
= \v_\lambda \* \ee{i_1} \jboxi{1} 
= \v_\lambda \* \jboxi{0} 
= \psi_1^{\lambda, \Lambda_{i_1}}(\v_{\Lambda_{i_1}})$.
So $\v_\lambda \* \jboxi{1}$ is \singular~ with $\xrightarrow{i_1}$
describing the only  \walks~ from the appropriate highest weight node
to it.

We compute, using the inductive hypothesis,
\begin{align*}
\psi_{k-m}^{\lambda, \mu} (\aA) 
&= \ff{i_1} (  \psi_1^{\lambda, \Lambda_{i_1}} \* {id}^{\*(k-1-m)})
\circ \psi_{k-1-m}^{\Lambda_{i_1}, \mu} (\ff{i_2} \cdots
\ff{i_k} \v_\mu) \\
&= \ff{i_1} (  \psi_1^{\lambda, \Lambda_{i_1}} \* {id}^{\*(k-1-m)})
(\v_{\Lambda_{i_1}} \* \jboxi{2} \*  \cdots \* \jboxi{k})
\\
&= \ff{i_1}  ( (\v_\lambda \* \jboxi{0}) \* \jboxi{2} \*  \cdots \* \jboxi{k})
\\
&= \v_\lambda \* \jboxi{1} \*  \cdots \* \jboxi{k}.
\end{align*}

3. This follows from 
Lemma \ref{lemma-eps-tensor}, and that $\psi_{k-m}^{\lambda,\mu}$ is
an embedding. 
Note that for $k \ge 3$, so long as $a_{i_1, i_2} < 0, a_{i_2, i_3} < 0$
the node will be \singular~ with \singular~ parent. 
\end{proof}

\section{Representation-theoretic interpretation in type $A$}
\label{sec-hecke}
In this paper, we studied a \singular~ node whose parent is also
\singular~ in highest weight crystals of finite and affine type.
The papers
\cite{Foda.LOTW.solvable}
\cite{FodaLOTW.unity}
characterized  all \singular~ nodes
in a level one highest weight crystal of type $A^{(1)}$
(and \cite{Foda.LOTW.Ariki.Koike} for higher levels)
 by their behavior under tensor
product of crystals, and they  gave a representation-theoretic interpretation
of these  \singular~ nodes  as answering  the Jantzen-Seitz problem.
These nodes correspond to irreducible modules of the
finite Hecke algebra $H_n^{\fin}$ of type $A$ that remain irreducible
on restriction from $H_n^{\fin}$ to $H_{n-1}^{\fin}$
(or for the Ariki-Koike (cyclotomic Hecke) algebras).
One may then ask: which
irreducible modules of 
$H_n^{\fin}$ remain irreducible
on restriction  to $H_{n-2}^{\fin}$?

The above Theorem \ref{thm-type} in type $A$ was motivated by the following
representation-theoretic fact, which addresses the question just posed. 
If an irreducible module $M$ of the {\it affine\/} Hecke algebra $H_n$ of type
$A$ is irreducible on restriction to $H_{n-2}$, then $M$ is one-dimensional
and either a trivial or Steinberg (sign) module. 
The main theorem of \cite{Groj.slp} says that
the above Hecke-theoretic statement
and
case \eqref{case-A}
of Theorem \ref{thm-type} 
are equivalent.
However, a purely representation-theoretic proof is as straightforward
as the crystal-theoretic  proofs above.

Compare the representation-theoretic translation  (given below)
of the crystal-theoretic
proof 
with the following direct proof communicated by Grojnowski.

Let $H_n$ denote the affine Hecke algebra of type $A$.
The algebra depends on a parameter $q$, and when we specialize
$q=1$, we recover   the  group algebra of the wreath product
of $\Z$ with 
the symmetric group.  We denote by $T_i$  the generator 
of $H_n$ that degenerates to the simple transposition $s_i = (i, i+1)$.

Let $M$ be an irreducible module of $H_n$,
and suppose $\Res_{H_{n-2}}^{H_n} M$ is an irreducible $H_{n-2}$-module.
In particular, the generator $T_{n-1}$ 
commutes with $H_{n-2}$ and so acts by a scalar on all of $M$,
where that scalar is $-1$ or $q$, as $(T_{n-1} + 1)(T_{n-1} - q) = 0$.
All of the $T_i$ are conjugate in $H_n$, so all the $T_i$ also act
by that same scalar on all of $M$.  In the case that scalar is $q$,
$M$ must have been a trivial module, and when it is $-1$
we have 
a Steinberg (sign)
module.  In particular,  $M$ is one-dimensional and  $\Res_{H_{k}}^{H_n} M$
is irreducible for all $k \le n$.
This argument is the correct explanation for the result,
but it is unclear what its interpretation
is in other types. 

\omitt{ 
Let $S_n$ denote the symmetric group.
We denote by $s_i$ the simple transposition $s_i = (i, i+1)$.
Let $M$ be an irreducible module of $S_n$,
and suppose $\Res_{S_{n-2}}^{S_n} M$ is an irreducible $S_{n-2}$-module.
As $s_{n-1}$ commutes with $S_{n-2}$, its eigenspaces give a decomposition
of $\Res_{S_{n-2}}^{S_n} M$.
As the restriction is irreducible, $s_{n-1}$ must act by the same
eigenvalue on all of $M$ and in particular acts as a scalar, either
$+1$ or $-1$.

All of the $s_i$ are conjugate in $S_n$, so all the $s_i$ also act
by that same scalar on all of $M$.  In the case that scalar is $1$ 
$M$ must have been a trivial module, and when it is $-1$
we have the sign module. 
In particular,  $M$ is one-dimensional and  $\Res_{S_{k}}^{S_n} M$
is irreducible for all $k \le n$.
This argument is the correct explanation for the result,
but it is unclear what its interpretation
is in other types. 
}

In contrast, here is the representation-theoretic version of the
crystal-theoretic proof  given in  case \eqref{case-A}
of Theorem \ref{thm-type}. 


We refer the reader to \cite{Groj.slp} and \cite{Groj.Vaz} for all
the definitions (as it is not the main focus of this paper).
Let $M$ be an irreducible module of $H_n$.
There are operators
$$ \eei : \Rep H_n \to \Rep H_{n-1}$$
that satisfy $\bigoplus_i \eei M = \Res_{H_{n-1}}^{H_n} M$.
Further, if $\eei^2 M = 0$ but $\eei M \neq 0$, then $\eei M$ is
an irreducible $H_{n-1}$-module, and conversely.
A node being singular corresponds to $\Res_{H_{n-1}}^{H_n} M$
being irreducible. 
Hence the hypotheses of  case \eqref{case-A} of Theorem \ref{thm-type}
correspond to the assumption that $\Res_{H_{n-1}}^{H_n} M = \eei M$
is irreducible, $\eee{i+1} \eei M$ is also irreducible, and
$\eej \eei M = 0$ for $j \neq i+1$. This implies $\Res_{H_{n-2}}^{H_n} M
= \eee{i+1} \eei M$
is irreducible.  We want to conclude that
$\eee{i+2} \eee{i+1} \eei M$ is also irreducible or zero (so
we need only show $\eee{i+2}^2 \eee{i+1} \eei M = 0$, and that
$\eej \eee{i+1} \eei M = 0$ for $j \neq i+2$).
These all follow from the fact, shown in  \cite{Groj.slp}, \cite{Groj.Vaz},
that the $\eei$ satisfy the Serre relations of type A.
(For ease of exposition, we omit the case where the parameter
$q$ appearing in the definition of $H_n$ is a second root of unity,
corresponding  to case \ref{case-Aone} of Theorem \ref{thm-type}.
We omitted this case in Grojnowski's direct proof above as well, where one
must confront the fact that the $T_i$ may not act semisimply.)
The proof here is very close to that of case \eqref{case-A} of Theorem
\ref{thm-type} and \eqref{case-1} of Corollary \ref{cor-serre},
as they both rely on the Serre relations.

We also point out that this statement is obvious for the
representation theory of the symmetric group in characteristic $0$.
Here, irreducible representations are indexed by partitions, and
the branching rule says the restriction of an irreducible can be
described by removing certain boxes from the partition. 
For the restriction of an irreducible module from $S_n$ to $S_{n-1}$
to be irreducible means its partition can have at most one removable
box, and hence be a rectangle.  But for that rectangle to share
the same property, the original shape must have been a single row or
column, hence our original representation was the trivial or sign module.
We remark that the combinatorics in prime characteristic are
appreciably different.

While for symmetric group modules in characteristic 0 this is
a classical fact, it seemed a surprising statement for
crystals: that two consecutive \singular~ nodes could determine
all of their ancestors, and that the perfect crystal
$B^{1,1}$ controls all the paths between that node and the
highest weight node.

\section{Exceptional types}
\label{sec-exceptional}
Corollaries \ref{cor-zero}, \ref{cor-parent}, \ref{cor-serre}
say that  in simply laced type, if 
$\varepsilon(\aA) = \Lambda_i, \varepsilon(\ei \aA) = \Lambda_j$,
and $\varepsilon_k(\ej \ei \aA) \neq 0$,
then we see
\begin{picture}(46,20)(-5,-10)
\multiput( 0,0)(20,0){3}{\circle{6}}
\multiput( 3,0)(20,0){2}{\line(1,0){14}}
\put(0,-5){\makebox(0,0)[t]{$i$}}
\put(20,-5){\makebox(0,0)[t]{$j$}}
\put(40,-5){\makebox(0,0)[t]{$k$}}
\end{picture}
in the Dynkin diagram. 
In particular $k \neq i$.
In classical types, the possible
$\v \xrightarrow{i_k} \cdots \xrightarrow{i_2}\xrightarrow{i_1} \aA$
are in correspondence with walks on $B^{1,1}$.
In exceptional types, one can also describe a directed graph
the corresponding walks must live on.  
The directed graph is dictated by equation \eqref{eq-weight}
and case 2 of Corollary \ref{cor-serre}.
They are very complicated to draw (planarly), so we only give pictures for $E_6$
below.  Just as in type $A$, the two graphs below differ by reversing
orientation of all arrows. 

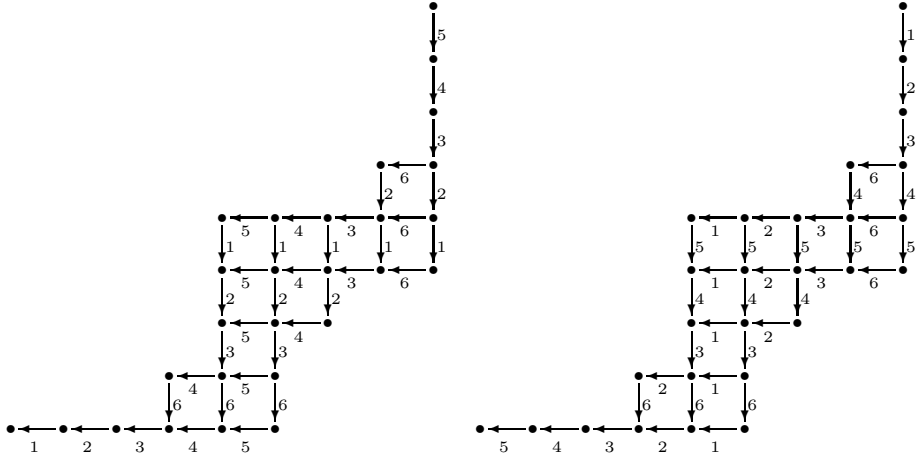
\begin{picture}(180,190)(25,0)
\multiput( 37,20)(20,0){5}{\vector(-1,0){14}}
\multiput( 20,20)(20,0){6}{\circle*{3}}
\multiput( 97,40)(20,0){2}{\vector(-1,0){14}}
\multiput( 80,40)(20,0){3}{\circle*{3}}
\multiput( 117,60)(20,0){2}{\vector(-1,0){14}}
\multiput( 100,60)(20,0){3}{\circle*{3}}
\multiput( 117,80)(20,0){4}{\vector(-1,0){14}}
\multiput( 100,80)(20,0){5}{\circle*{3}}
\multiput( 117,100)(20,0){4}{\vector(-1,0){14}}
\multiput( 100,100)(20,0){5}{\circle*{3}}
\put( 177,120){\vector(-1,0){14}}
\multiput( 160,120)(20,0){2}{\circle*{3}}
\multiput( 180,177)(0,-20){5}{\vector(0,-1){14}}
\multiput( 180,180)(0,-20){3}{\circle*{3}}
\multiput( 160,117)(0,-20){2}{\vector(0,-1){14}}
\multiput( 140,97)(0,-20){2}{\vector(0,-1){14}}
\multiput( 120,97)(0,-20){4}{\vector(0,-1){14}}
\multiput( 100,97)(0,-20){4}{\vector(0,-1){14}}
\put( 80,37){\vector(0,-1){14}}
\put(29,15){\makebox(0,0)[t]{${\scriptscriptstyle 1}$}}
\put(49,15){\makebox(0,0)[t]{${\scriptscriptstyle 2}$}}
\put(69,15){\makebox(0,0)[t]{${\scriptscriptstyle 3}$}}
\put(89,15){\makebox(0,0)[t]{ ${\scriptscriptstyle 4}$ }}
\put(109,15){\makebox(0,0)[t]{${\scriptscriptstyle 5}$}}
\put(89,37){\makebox(0,0)[t]{ ${\scriptscriptstyle 4}$ }}
\put(109,37){\makebox(0,0)[t]{${\scriptscriptstyle 5}$}}
\put(109,57){\makebox(0,0)[t]{${\scriptscriptstyle 5}$}}
\put(129,57){\makebox(0,0)[t]{${\scriptscriptstyle 4}$}}
\put(109,77){\makebox(0,0)[t]{${\scriptscriptstyle 5}$}}
\put(129,77){\makebox(0,0)[t]{${\scriptscriptstyle 4}$}}
\put(149,77){\makebox(0,0)[t]{${\scriptscriptstyle 3}$}}
\put(169,77){\makebox(0,0)[t]{${\scriptscriptstyle 6}$}}
\put(109,97){\makebox(0,0)[t]{${\scriptscriptstyle 5}$}}
\put(129,97){\makebox(0,0)[t]{${\scriptscriptstyle 4}$}}
\put(149,97){\makebox(0,0)[t]{${\scriptscriptstyle 3}$}}
\put(169,97){\makebox(0,0)[t]{${\scriptscriptstyle 6}$}}
\put(169,117){\makebox(0,0)[t]{${\scriptscriptstyle 6}$}}
\multiput( 83,31)(20,0){3}{\makebox(0,0)[t]{${\scriptscriptstyle 6}$}}
\multiput( 103,51)(20,0){2}{\makebox(0,0)[t]{${\scriptscriptstyle 3}$}}
\multiput( 103,71)(20,0){3}{\makebox(0,0)[t]{${\scriptscriptstyle 2}$}}
\multiput( 103,91)(20,0){5}{\makebox(0,0)[t]{${\scriptscriptstyle 1}$}}
\multiput( 163,111)(20,0){2}{\makebox(0,0)[t]{${\scriptscriptstyle 2}$}}
\put(183,131){\makebox(0,0)[t]{${\scriptscriptstyle 3}$}}
\put(183,151){\makebox(0,0)[t]{${\scriptscriptstyle 4}$}}
\put(183,171){\makebox(0,0)[t]{${\scriptscriptstyle 5}$}}
\end{picture}
\begin{picture}(180,190)(30,0)
\multiput( 37,20)(20,0){5}{\vector(-1,0){14}}
\multiput( 20,20)(20,0){6}{\circle*{3}}
\multiput( 97,40)(20,0){2}{\vector(-1,0){14}}
\multiput( 80,40)(20,0){3}{\circle*{3}}
\multiput( 117,60)(20,0){2}{\vector(-1,0){14}}
\multiput( 100,60)(20,0){3}{\circle*{3}}
\multiput( 117,80)(20,0){4}{\vector(-1,0){14}}
\multiput( 100,80)(20,0){5}{\circle*{3}}
\multiput( 117,100)(20,0){4}{\vector(-1,0){14}}
\multiput( 100,100)(20,0){5}{\circle*{3}}
\put( 177,120){\vector(-1,0){14}}
\multiput( 160,120)(20,0){2}{\circle*{3}}
\multiput( 180,177)(0,-20){5}{\vector(0,-1){14}}
\multiput( 180,180)(0,-20){3}{\circle*{3}}
\multiput( 160,117)(0,-20){2}{\vector(0,-1){14}}
\multiput( 140,97)(0,-20){2}{\vector(0,-1){14}}
\multiput( 120,97)(0,-20){4}{\vector(0,-1){14}}
\multiput( 100,97)(0,-20){4}{\vector(0,-1){14}}
\put( 80,37){\vector(0,-1){14}}
\put(29,15){\makebox(0,0)[t]{${\scriptscriptstyle 5}$}}
\put(49,15){\makebox(0,0)[t]{${\scriptscriptstyle 4}$}}
\put(69,15){\makebox(0,0)[t]{${\scriptscriptstyle 3}$}}
\put(89,15){\makebox(0,0)[t]{ ${\scriptscriptstyle 2}$ }}
\put(109,15){\makebox(0,0)[t]{${\scriptscriptstyle 1}$}}
\multiput( 109,37)(0,20){4}{\makebox(0,0)[t]{${\scriptscriptstyle 1}$}}
\put(89,37){\makebox(0,0)[t]{ ${\scriptscriptstyle 2}$ }}
\multiput( 129,57)(0,20){3}{\makebox(0,0)[t]{${\scriptscriptstyle 2}$}}
\put(149,77){\makebox(0,0)[t]{${\scriptscriptstyle 3}$}}
\multiput( 169,77)(0,20){3}{\makebox(0,0)[t]{${\scriptscriptstyle 6}$}}
\put(149,97){\makebox(0,0)[t]{${\scriptscriptstyle 3}$}}
\multiput( 83,31)(20,0){3}{\makebox(0,0)[t]{${\scriptscriptstyle 6}$}}
\multiput( 103,51)(20,0){2}{\makebox(0,0)[t]{${\scriptscriptstyle 3}$}}
\multiput( 103,71)(20,0){3}{\makebox(0,0)[t]{${\scriptscriptstyle 4}$}}
\multiput( 103,91)(20,0){5}{\makebox(0,0)[t]{${\scriptscriptstyle 5}$}}
\multiput( 163,111)(20,0){2}{\makebox(0,0)[t]{${\scriptscriptstyle 4}$}}
\put(183,131){\makebox(0,0)[t]{${\scriptscriptstyle 3}$}}
\put(183,151){\makebox(0,0)[t]{${\scriptscriptstyle 2}$}}
\put(183,171){\makebox(0,0)[t]{${\scriptscriptstyle 1}$}}
\end{picture}

\section{Appendix}
\label{sec-appendix}

The crystal graphs $B^{1,1}$ are listed  below.

We also need $B^{n-1,1}$ in type $A$:
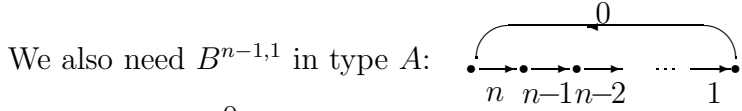
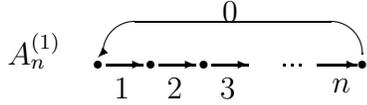
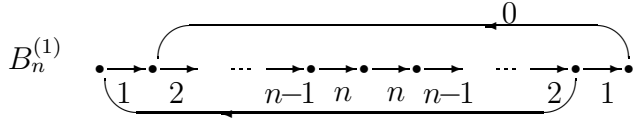
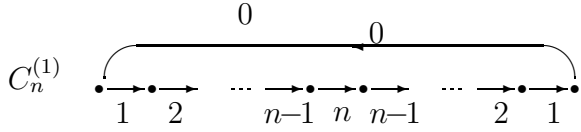
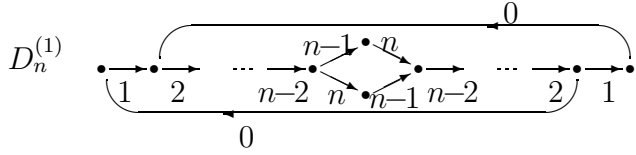
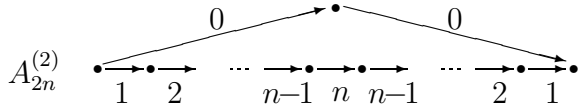
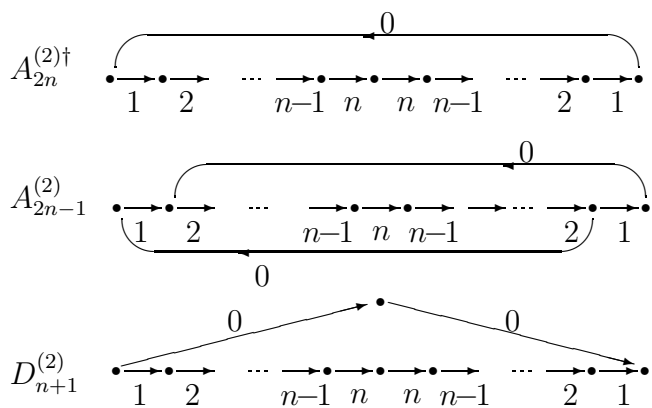
\begin{picture}(150,38)(-10,0)
\multiput( 0,0)(20,0){3}{\circle*{3}}
\multiput( 3,0)(20,0){3}{\vector(1,0){14}}
\multiput(70,0)(3,0){3}{\line(1,0){1}} 
\put(100,0){\circle*{3}}
\put(83,0){\vector(1,0){15}}
\put(51,4){\oval(98,25)[t]}
\put(47,16){\vector(-1,0){3}}
\put(9,-7){\makebox(0,0)[t]{$n$}}
\put(29,-5){\makebox(0,0)[t]{$n\! \!-\! \!1$}}
\put(49,-5){\makebox(0,0)[t]{$n\! \!-\! \!2$}}
\put(92,-5){\makebox(0,0)[t]{$1$}}
\put(50,25){\makebox(0,0)[t]{$0$}}
\end{picture}

$A_n^{(1)}$
\begin{picture}(150,35)(-10,0)
\multiput( 0,0)(20,0){3}{\circle*{3}}
\multiput( 3,0)(20,0){3}{\vector(1,0){14}}
\multiput(70,0)(3,0){3}{\line(1,0){1}} 
\put(100,0){\circle*{3}}
\put(83,0){\vector(1,0){15}}
\put(51,4){\oval(98,25)[t]}
\put(3,6){\vector(-1,-1){3}}
\put(9,-5){\makebox(0,0)[t]{$1$}}
\put(29,-5){\makebox(0,0)[t]{$2$}}
\put(49,-5){\makebox(0,0)[t]{$3$}}
\put(92,-5){\makebox(0,0)[t]{$n$}}
\put(50,24){\makebox(0,0)[t]{$0$}}
\end{picture}
%

%
$B_n^{(1)}$
\begin{picture}(180,40)(-10,0)
\multiput( 0,0)(20,0){2}{\circle*{3}}
\multiput(80,0)(20,0){3}{\circle*{3}}
\multiput(180,0)(20,0){2}{\circle*{3}}
\multiput( 3,0)(20,0){2}{\vector(1,0){14}}
\multiput(63,0)(20,0){4}{\vector(1,0){14}}
\multiput(163,0)(20,0){2}{\vector(1,0){14}}
\multiput(50,0)(3,0){3}{\line(1,0){1}} 
\multiput(150,0)(3,0){3}{\line(1,0){1}} 
\put(91,-4){\oval(178,25)[b]}
\put(111,4){\oval(178,25)[t]}
\put(149,16){\vector(-1,0){3}}
\put(49,-17){\vector(-1,0){3}}
\put(9,-5){\makebox(0,0)[t]{$1$}}
\put(29,-5){\makebox(0,0)[t]{$2$}}
\put(72,-5){\makebox(0,0)[t]{$n\! \!-\! \!1$}}
\put(92,-7){\makebox(0,0)[t]{$n$}}
\put(112,-7){\makebox(0,0)[t]{$n$}}
\put(132,-5){\makebox(0,0)[t]{$n\! \!-\! \!1$}}
\put(172,-5){\makebox(0,0)[t]{$2$}}
\put(192,-5){\makebox(0,0)[t]{$1$}}
\put(155,25){\makebox(0,0)[t]{$0$}}
\put(55,-22){\makebox(0,0)[t]{$0$}}
\end{picture}
%

%
$C_n^{(1)}$
\begin{picture}(200,52)(-10,0)
\multiput( 0,0)(20,0){2}{\circle*{3}}
\multiput(80,0)(20,0){2}{\circle*{3}}
\multiput(160,0)(20,0){2}{\circle*{3}}
\multiput( 3,0)(20,0){2}{\vector(1,0){14}}
\multiput(63,0)(20,0){3}{\vector(1,0){14}}
\multiput(143,0)(20,0){2}{\vector(1,0){14}}
\multiput(50,0)(3,0){3}{\line(1,0){1}} 
\multiput(130,0)(3,0){3}{\line(1,0){1}} 
\put(91,4){\oval(178,25)[t]}
\put(99,16){\vector(-1,0){3}}
\put(9,-5){\makebox(0,0)[t]{$1$}}
\put(29,-5){\makebox(0,0)[t]{$2$}}
\put(72,-5){\makebox(0,0)[t]{$n\! \!-\! \!1$}}
\put(92,-7){\makebox(0,0)[t]{$n$}}
\put(112,-5){\makebox(0,0)[t]{$n\! \!-\! \!1$}}
\put(152,-5){\makebox(0,0)[t]{$2$}}
\put(172,-5){\makebox(0,0)[t]{$1$}}
\put(105,25){\makebox(0,0)[t]{$0$}}
\end{picture}

%
$D_n^{(1)}$
\begin{picture}(215,45)(-10,0)
\multiput( 0,0)(20,0){2}{\circle*{3}}
\multiput(80,0)(40,0){2}{\circle*{3}}
\multiput(100,-10)(0,20){2}{\circle*{3}}
\multiput(180,0)(20,0){2}{\circle*{3}}
\multiput( 3,0)(20,0){2}{\vector(1,0){14}}
\multiput(63,0)(60,0){2}{\vector(1,0){14}}
\multiput(163,0)(20,0){2}{\vector(1,0){14}}
\multiput(50,0)(3,0){3}{\line(1,0){1}} 
\multiput(150,0)(3,0){3}{\line(1,0){1}} 
\put(91,-4){\oval(178,25)[b]}
\put(111,4){\oval(178,25)[t]}
\put(149,16){\vector(-1,0){3}}
\put(49,-17){\vector(-1,0){3}}
\put(9,-5){\makebox(0,0)[t]{$1$}}
\put(29,-5){\makebox(0,0)[t]{$2$}}
\put(69,-5){\makebox(0,0)[t]{$n\! \!-\! \!2$}}
\put(133,-5){\makebox(0,0)[t]{$n\! \!-\! \!2$}}
\put(83,1){\vector(2,1){14}}
\put(83,-1){\vector(2,-1){14}}
\put(103,9){\vector(2,-1){14}}
\put(103,-9){\vector(2,1){14}}
\put(172,-5){\makebox(0,0)[t]{$2$}}
\put(192,-5){\makebox(0,0)[t]{$1$}}
\put(155,25){\makebox(0,0)[t]{$0$}}
\put(55,-22){\makebox(0,0)[t]{$0$}}
\put(86,13){\makebox(0,0)[t]{$n\! \!-\! \!1$}}
\put(89,-8){\makebox(0,0)[t]{$n$}}
\put(111,-7){\makebox(0,0)[t]{$n\! \!-\! \!1$}}
\put(109,12){\makebox(0,0)[t]{$n$}}
\end{picture}

%
$A_{2n}^{(2)}$
\begin{picture}(200,65)(-10,-5)
\multiput( 0,0)(20,0){2}{\circle*{3}}
\multiput(80,0)(20,0){2}{\circle*{3}}
\multiput(160,0)(20,0){2}{\circle*{3}}
\multiput( 3,0)(20,0){2}{\vector(1,0){14}}
\multiput(63,0)(20,0){3}{\vector(1,0){14}}
\multiput(143,0)(20,0){2}{\vector(1,0){14}}
\multiput(50,0)(3,0){3}{\line(1,0){1}} 
\multiput(130,0)(3,0){3}{\line(1,0){1}} 
\put(90,23){\circle*{3}}
\put(2,2){\vector(4,1){84}}
\put(93,24){\vector(4,-1){84}}
\put(9,-5){\makebox(0,0)[t]{$1$}}
\put(29,-5){\makebox(0,0)[t]{$2$}}
\put(72,-5){\makebox(0,0)[t]{$n\! \!-\! \!1$}}
\put(92,-7){\makebox(0,0)[t]{$n$}}
\put(112,-5){\makebox(0,0)[t]{$n\! \!-\! \!1$}}
\put(152,-5){\makebox(0,0)[t]{$2$}}
\put(172,-5){\makebox(0,0)[t]{$1$}}
\put(135,23){\makebox(0,0)[t]{$0$}}
\put(45,23){\makebox(0,0)[t]{$0$}}
\end{picture}

%
%
$A_{2n}^{(2)\dagger}$
\begin{picture}(200,45)(-10,0)
\multiput( 0,0)(20,0){2}{\circle*{3}}
\multiput(80,0)(20,0){3}{\circle*{3}}
\multiput(180,0)(20,0){2}{\circle*{3}}
\multiput( 3,0)(20,0){2}{\vector(1,0){14}}
\multiput(63,0)(20,0){4}{\vector(1,0){14}}
\multiput(163,0)(20,0){2}{\vector(1,0){14}}
\multiput(50,0)(3,0){3}{\line(1,0){1}} 
\multiput(150,0)(3,0){3}{\line(1,0){1}} 
\put(101,4){\oval(198,25)[t]}
\put(99,16){\vector(-1,0){3}}
\put(9,-5){\makebox(0,0)[t]{$1$}}
\put(29,-5){\makebox(0,0)[t]{$2$}}
\put(72,-5){\makebox(0,0)[t]{$n\! \!-\! \!1$}}
\put(92,-7){\makebox(0,0)[t]{$n$}}
\put(112,-7){\makebox(0,0)[t]{$n$}}
\put(132,-5){\makebox(0,0)[t]{$n\! \!-\! \!1$}}
\put(172,-5){\makebox(0,0)[t]{$2$}}
\put(192,-5){\makebox(0,0)[t]{$1$}}
\put(105,25){\makebox(0,0)[t]{$0$}}
\end{picture}

%
$A_{2n-1}^{(2)}$
\begin{picture}(210,45)(-7,0)
\multiput( 0,0)(20,0){2}{\circle*{3}}
\multiput(90,0)(20,0){2}{\circle*{3}}
\multiput(180,0)(20,0){2}{\circle*{3}}
\multiput( 3,0)(20,0){2}{\vector(1,0){14}}
\multiput(73,0)(20,0){4}{\vector(1,0){14}}
\multiput(163,0)(20,0){2}{\vector(1,0){14}}
\multiput(50,0)(3,0){3}{\line(1,0){1}} 
\multiput(150,0)(3,0){3}{\line(1,0){1}} 
\put(91,-4){\oval(178,25)[b]}
\put(111,4){\oval(178,25)[t]}
\put(149,16){\vector(-1,0){3}}
\put(49,-17){\vector(-1,0){3}}
\put(9,-5){\makebox(0,0)[t]{$1$}}
\put(29,-5){\makebox(0,0)[t]{$2$}}
\put(80,-5){\makebox(0,0)[t]{$n\! \!-\! \!1$}}
\put(100,-7){\makebox(0,0)[t]{$n$}}
\put(120,-5){\makebox(0,0)[t]{$n\! \!-\! \!1$}}
\put(172,-5){\makebox(0,0)[t]{$2$}}
\put(192,-5){\makebox(0,0)[t]{$1$}}
\put(155,25){\makebox(0,0)[t]{$0$}}
\put(55,-22){\makebox(0,0)[t]{$0$}}
\end{picture}

%
$D_{n+1}^{(2)}$
\begin{picture}(190,62)(-10,-5)
\multiput( 0,0)(20,0){2}{\circle*{3}}
\multiput(80,0)(20,0){3}{\circle*{3}}
\multiput(180,0)(20,0){2}{\circle*{3}}
\multiput( 3,0)(20,0){2}{\vector(1,0){14}}
\multiput(63,0)(20,0){4}{\vector(1,0){14}}
\multiput(163,0)(20,0){2}{\vector(1,0){14}}
\multiput(50,0)(3,0){3}{\line(1,0){1}} 
\multiput(150,0)(3,0){3}{\line(1,0){1}} 
\put(9,-5){\makebox(0,0)[t]{$1$}}
\put(29,-5){\makebox(0,0)[t]{$2$}}
\put(72,-5){\makebox(0,0)[t]{$n\! \!-\! \!1$}}
\put(92,-7){\makebox(0,0)[t]{$n$}}
\put(112,-7){\makebox(0,0)[t]{$n$}}
\put(132,-5){\makebox(0,0)[t]{$n\! \!-\! \!1$}}
\put(172,-5){\makebox(0,0)[t]{$2$}}
\put(192,-5){\makebox(0,0)[t]{$1$}}
\put(100,26){\circle*{3}}
\put(2,2){\vector(4,1){93}}
\put(103,26){\vector(4,-1){93}}
\put(150,23){\makebox(0,0)[t]{$0$}}
\put(45,23){\makebox(0,0)[t]{$0$}}
\end{picture}


\bibliographystyle{alpha}

\begin{thebibliography}{FLO{\etalchar{+}}98b}

\bibitem[FLO{\etalchar{+}}98a]{Foda.LOTW.solvable}
O.~Foda, B.~Leclerc, M.~Okado, J.-Y. Thibon, and T.~A. Welsh.
\newblock Combinatorics of solvable lattice models, and modular representations
  of {H}ecke algebras.
\newblock In {\em Geometric analysis and Lie theory in mathematics and
  physics}, volume~11 of {\em Austral. Math. Soc. Lect. Ser.}, pages 243--290.
  Cambridge Univ. Press, Cambridge, 1998.

\bibitem[FLO{\etalchar{+}}98b]{FodaLOTW.unity}
Omar Foda, Bernard Leclerc, Masato Okado, Jean-Yves Thibon, and Trevor~A.
  Welsh.
\newblock R{SOS} models and {J}antzen-{S}eitz representations of {H}ecke
  algebras at roots of unity.
\newblock {\em Lett. Math. Phys.}, 43(1):31--42, 1998.

\bibitem[FLO{\etalchar{+}}99]{Foda.LOTW.Ariki.Koike}
Omar Foda, Bernard Leclerc, Masato Okado, Jean-Yves Thibon, and Trevor~A.
  Welsh.
\newblock Branching functions of {$A\sp {(1)}\sb {n-1}$} and {J}antzen-{S}eitz
  problem for {A}riki-{K}oike algebras.
\newblock {\em Adv. Math.}, 141(2):322--365, 1999.

\bibitem[Gro]{Groj.slp}
I.~Grojnowski.
\newblock {Affine ${\mathfrak sl}_p$ controls the representation theory of the
  symmetric group and related Hecke algebras}.

\bibitem[GV01]{Groj.Vaz}
I.~Grojnowski and M.~Vazirani.
\newblock Strong multiplicity one theorems for affine {H}ecke algebras of type
  {A}.
\newblock {\em Transform. Groups}, 6(2):143--155, 2001.

\bibitem[Kas93]{Kashiwara.global}
Masaki Kashiwara.
\newblock Global crystal bases of quantum groups.
\newblock {\em Duke Math. J.}, 69(2):455--485, 1993.

\bibitem[Kas95]{Kashiwara.on.crystal.bases}
Masaki Kashiwara.
\newblock On crystal bases.
\newblock In {\em Representations of groups (Banff, AB, 1994)}, volume~16 of
  {\em CMS Conf. Proc.}, pages 155--197. Amer. Math. Soc., Providence, RI,
  1995.

\bibitem[KK98]{KK.core}
Seok-Jin Kang and Masaki Kashiwara.
\newblock Quantized affine algebras and crystals with core.
\newblock {\em Comm. Math. Phys.}, 195(3):725--740, 1998.

\bibitem[KKM{\etalchar{+}}92]{KKMMNN.perfect}
Seok-Jin Kang, Masaki Kashiwara, Kailash~C. Misra, Tetsuji Miwa, Toshiki
  Nakashima, and Atsushi Nakayashiki.
\newblock Perfect crystals of quantum affine {L}ie algebras.
\newblock {\em Duke Math. J.}, 68(3):499--607, 1992.

\bibitem[Ste03]{Stembridge.crystal}
John~R. Stembridge.
\newblock A local characterization of simply-laced crystals.
\newblock {\em Trans. Amer. Math. Soc.}, 355(12):4807--4823 (electronic), 2003.

\end{thebibliography}

\end{document}